\documentclass[12pt,reqno]{amsart}
\usepackage{xypic, amssymb, amsfonts, amsbsy, amsthm, amsmath, amscd, latexsym, stmaryrd, epic, eepic,eucal}

\DeclareMathAlphabet{\mathpzc}{OT1}{pzc}{m}{it}

\newenvironment{dem}{\begin{proof}[\bf Proof]}{\end{proof}}

\newtheorem{theorem}{\bf Theorem}[section]
\newtheorem{lemma}[theorem]{\bf Lemma}
\newtheorem{propos}[theorem]{\bf Proposition}
\newtheorem{corol}[theorem]{\bf Corollary}
\newtheorem{claim}[theorem]{\bf Claim}

\theoremstyle{definition}
\newtheorem{defi}[theorem]{\bf Definition}
\newtheorem{oss}[theorem]{\bf Remark}

\newtheorem{exm}[theorem]{\bf Example}

\newcommand{\A}{\mathbb A}
\newcommand{\B}{\text{B}}

\newcommand{\C}{\mathbb C}
\newcommand{\Cg}{\mathfrak C}
\newcommand{\F}{\mathcal F}
\newcommand{\Gm}{\mathbb G_{\textbf{m}}}
\newcommand{\Ga}{\mathbb G_{\textbf{a}}}

\newcommand{\M}{\mathfrak M}

\newcommand{\Ms}{\mathcal M}
\newcommand{\N}{\mathbb N}
\newcommand{\Nf}{\mathcal N}
\newcommand{\Of}{\mathcal O}
\newcommand{\Pro}{\mathbb P}
\newcommand{\Q}{\mathbb Q}

\newcommand{\St}{\mathcal S}

\newcommand{\Z}{\mathbb Z}
\newcommand{\aut}{\text{Aut}}

\newcommand{\az}{\gamma}

\newcommand{\cart}{\ar @{} [dr] |{\Box}}
\newcommand{\cat}{\mathpzc{Sch}_{\mathbb C}}

\newcommand{\cm}{\mathpzc c}
\newcommand{\cu}{\mathpzc C}

\newcommand{\hol}{\frac{\partial}{\partial z}}
\newcommand{\id}{\text{id}}

\newcommand{\km}{\mathpzc k}

\newcommand{\nm}{\mathpzc n}

\newcommand{\rat}{\otimes \mathbb Q}

\newcommand{\sla}{\mathfrak{sl}}
\newcommand{\spc}{\text{Spec} \mathbb C}
\newcommand{\spe}{\text{Spec}}
\newcommand{\spk}{\text{Spec}\Omega}

\newcommand{\unoa}{\xymatrix { *=0{\bullet} \ar@{-}[r] & *=0{\bullet} }}
\newcommand{\duea}{\xymatrix { *=0{\bullet} \ar@{-}[r] & *=0{\bullet} \ar@{-}[r] &*=0{\bullet} }}

\newcommand{\trea}{\xymatrix { *=0{\bullet} \ar@{-}[r] & *=0{\bullet} \ar@{-}[r] &*=0{\bullet} \ar@{-}[r] &*=0{\bullet} }}

\newcommand{\treb} {\xymatrix @R=8pt { & *=0{\bullet} \ar@{-}[dd] & \\ & &  \\ & *=0{\bullet} \ar@{-}[dl] \ar@{-}[dr] & \\ *=0{\bullet} & & *=0{\bullet}}}

\newcommand{\trebg}{\xymatrix @R=8pt { & *=0{} \ar@{-}[4,0] & \\ *=0{} \ar@{-}[0,2] & & *=0{} \\ *=0{} \ar@{-}[0,2] & & *=0{}\\ *=0{} \ar@{-}[0,2] & & *=0{}\\ & *=0{} &}}

\newcommand{\quattroa}{\xymatrix { *=0{\bullet} \ar@{-}[r] & *=0{\bullet} \ar@{-}[r] &*=0{\bullet} \ar@{-}[r] &*=0{\bullet} \ar@{-}[r] &*=0{\bullet} }}

\newcommand{\quattrob}{\xymatrix {& *=0{\bullet}  \ar@{-}[d]& \\ *=0{\bullet} \ar@{-}[r] & *=0{\bullet} \ar@{-}[r] \ar@{-}[d] & *=0{\bullet}\\ & *=0{\bullet} &}}

\newcommand{\quattroc}{\xymatrix @R=4pt {*=0{\bullet} \ar@{-}[dr] & & &\\ & *=0{\bullet} \ar@{-}[r] &*=0{\bullet} \ar@{-}[r] &*=0{\bullet}\\ *=0{\bullet} \ar@{-}[ur] & & &}}

\newcommand{\quattrocg}{\xymatrix @R=8pt { & *=0{} \ar@{-}[4,0] & & &\\ *=0{} \ar@{-}[0,2] & & *=0{} &*=0{} \ar@{-}[2,0] & \\ *=0{} \ar@{-}[0,4] & & & &*=0{}\\ *=0{} \ar@{-}[0,2] & & *=0{} &*=0{} &\\ & *=0{} & & &}}

\newcommand{\cinquenz}{\xymatrix { *=0{\bullet} \ar@{-}[d] & *=0{\bullet} \ar@{-}[d] \\ *=0{\bullet} \ar@{-}[r] & *=0{\bullet} \\
*=0{\bullet} \ar@{-}[u] & *=0{\bullet} \ar@{-}[u]}}

\input xy
\xyoption{all}

\begin{document}

\parindent = 0 pt

\title[Chow Ring of Rational Curves]{ON THE CHOW RING OF THE STACK OF RATIONAL NODAL CURVES}
\author{Damiano Fulghesu}
\address{Department of Mathematics,
University of Missouri,
Columbia, MO 65211}
\email{damiano@math.missouri.edu} 

\maketitle

\begin{abstract}
The goal of this paper is to compute the rational Chow ring of the stack $\M_{0}^{\leq 3}$ consisting of nodal curves of genus $0$ with at most $3$ nodes: it is a $\Q$-algebra with $10$ generators and $11$ relations.
\end{abstract}

\section{Introduction}

\medskip

\subsection{General Background}

Intersection theory on moduli spaces of stable curves started at the beginning of the '80's with Mumford's paper \cite{mumenum}, where he laid the foundations and carried out the first calculations. Many people have contributed to the theory after this (such as Faber \cite{Fab}, Witten and Kontsevich), building an imposing structure.

The foundations of intersection theory on Deligne-Mumford stacks have been developed by Gillet \cite{Gil} and Vistoli \cite{Vis}. The first step towards an intersection theory on general Artin stacks (like those that arise from looking at unstable curves) was the equivariant intersection theory that Edidin and Graham \cite{EGRR} developed, following an idea of Totaro \cite{tot}. Their theory associates a commutative graded Chow ring $A^*(\mathcal M)$ with every smooth quotient stack $\mathcal M$ of finite type over a field.

Unfortunately many stacks of geometric interest are not known to be quotient stacks (the general question of when a stack is a quotient stack is not completely understood; anyway there is a first answer in \cite{EHKV}). Later, A. Kresch \cite{kre} developed an intersection theory for general Artin stacks; in particular, he associates a Chow ring $A^*(\mathcal M)$ with every smooth Artin stack $\mathcal M$ locally of finite type over a field, provided some technical conditions hold, which are satisfied in particular for stacks of pointed nodal curves of fixed genus.

Since there is not yet a theory of Chow rings of such stacks that extends the theory of stacks of stable curves, there does not seem to be much else to do than look at specific examples. The first example is the stack $\M_{0}$ of nodal connected curves of genus $0$. However, even this case turns out  to be extremely complicated.

\medskip

\subsection{Description of contents}

We start with giving a brief description of the Artin stack  $\M_{0}$ of nodal connected curves of genus $0$. In (\ref{HV}) we give an example of a smooth projective surface $S$ with a family $C \to S$ of nodal curves of genus $0$ with at most $2$ nodes in each fiber, in which $C$ is an algebraic space not scheme-like. 

Then we introduce one of the basic tools: the stratification by nodes. Denote by $\M_{0}^{n}$ the smooth locally closed substack of $\M_{0}$ consisting of curves with exactly $n$ nodes. For each $n$, the closed substack $\M_{0}^{\geq n}:= \bigcup_{i \geq n}\M_{0}^{i}$ (consisting of families of curves whose geometric fibers have at least $n$ nodes) admits a regular embedding into $\M_0$ and it has codimension $n$ (see Proposition \ref{stratreg}). This gives a stratification of  $\M_0$. Each stratum has a decomposition into irreducible substacks that classify the topological type of a rational nodal curve. 

We can represent this topological type by a tree, with the  vertices corresponding to components and edges corresponding to nodes. Therefore for each tree $\Gamma$, there is a natural smooth locally closed substack $\M_{0}^{\Gamma}$ defined as the category of families of curves whose fibers have topological type corresponding to $\Gamma$.

In Section \ref{strata} we give an explicit description of stacks $\M^{\Gamma}_0$ for particular $\Gamma$: if $\Gamma$ comes from a curve whose components have no more that 3 nodes, then we have $\M^{\Gamma}_0 = \B \aut (C)$.

In particular we can stratify $\M^{\leq 3}_0$ with smooth quotient stacks. This allows to apply Kresch's theory and define a Chow ring $A^*(\M^{\leq 3}_{0})$. 

Our technique is to compute $A^*(\M_{0}^{\leq n})\rat $ for $n \leq 3$ by induction on $n$. For $n = 0$ we have $\M_{0}^{0} = \B \Pro GL_{2}$, and this case is well understood.

The inductive step is based on the following fact: if $n \leq 4$, then the top Chern class of the normal bundle of $\M_{0}^{n}$ into $\M_{0}^{\leq n}$ is not a $0$-divisor in $A^*(\M_{0}^{n})\rat$. As a consequence, by an elementary algebraic Lemma (\ref{incollamento}) we can reconstruct the ring $A^*(\M_{0}^{\leq n})\rat $ from the rings $A^*(\M_{0}^{\leq n-1})\rat$ and $A^*(\M_{0}^{n})\rat$ (Proposition \ref{ringstratum}), provided that we have an explicit way of extending each class in $A^*(\M_{0}^{\leq n-1})\rat$ to a class in $A^*(\M_{0}^{\leq n})\rat$ (Section \ref{restrcla}), and then computing the restriction of this extension to $A^*(\M_{0}^{n})\rat$.

We are able to describe completely two kinds of classes in $A^*(\M^{\leq 3}_0)$:

{\bf Strata classes.} For each tree $\Gamma$ with $\delta$ edges and maximal multiplicity 3, we get a class $\gamma_{\Gamma}$ of codimension $\delta$ in $\M_{0}$: the class of the closure of $\M_{0}^{\Gamma}$ in $\M_{0}$. We also compute the restriction of each $\gamma_{\Gamma}$ to all the rings $A^*(\M_{0}^{\Gamma'})\rat$ for any other tree $\Gamma'$ with maximal multiplicity 3 (Proposition \ref{refrestr}).

{\bf Mumford classes.} These should be defined as follows (Section \ref{mumcla}): let $\cu \xrightarrow{\Pi} \M_{0}$ be the universal curve. Call $K \in A^1(\cu)\rat$ the first Chern class of the relative dualizing sheaf $\omega_{\cu/ \M_{0}}$, and set
   $$
   \km_{i} = \Pi_{*}(K^{i+1}).
   $$
However, here we encounter an unfortunate technical problem: the morphism $\cu \to \M_{0}$ is not projective (not even represented by schemes, as Example \ref{HV} shows) and in Kresch's theory one does not have arbitrary proper pushforwards, only pushforwards along projective morphisms.

We are able to circumvent this problem only for curves with at most $3$ nodes. This works as follows. We show in Proposition (\ref{dualsh}) that the pushforward $\Pi_{*}\omega^{\vee}_{\cu / \M_{0}}$ is a locally free sheaf of rank $3$ on the the open substack $\M_{0}^{\leq 3}$ (this fails for curves with $4$ nodes). Hence we have Chern classes  $\cm_{1}$, $\cm_{2}$ and $\cm_{3}$ of $\Pi_{*}\omega^{\vee}_{\cu / \M_{0}}$ in $A^*(\M_{0}^{\leq 3})$.

We have a pushforward $\Pi_{*} A^*(\mathcal C^{\Gamma})\rat \to A^*(\M_{0}^{\Gamma})\rat$ along the restriction $\mathcal C^{\Gamma} \xrightarrow{\Pi} \M_{0}^{\Gamma}$ of the universal curve, because the stacks involved are quotient stacks, and arbitrary proper pushforwards exists in the theory of Edidin and Graham. This allows, with Grothendieck-Riemann-Roch, to compute the restriction of the Mumford classes to each $A^*(\M_{0}^{\Gamma})\rat$ (Proposition \ref{defmumstrat}), even without knowing that the Mumford classes exist. It turns out that for each tree $\Gamma$ with at most $3$ nodes, the Mumford classes in $A^*(\M_{0}^{\Gamma})\rat$ are polynomials in the restrictions of $\cm_{1}$, $\cm_{2}$ and $\cm_{3}$, thought of as elementary symmetric polynomials in three variables. Then we define the classes $\km_{i} \in A^*(\M_{0}^{\leq 3})\rat$ as the suitable polynomials in the $\cm_{1}$, $\cm_{2}$ and $\cm_{3}$ (Definition \ref{defmumpol}).

Since $A^*(\M_{0}^{\leq3})\rat$ injects into the product of $A^*(\M_{0}^{\Gamma})\rat$ over trees with at most three nodes (Proposition \ref{injleqtre}), this gives the right definition.

In the last Section we put everything together, and we calculate $A^*(\M_{0}^{\leq 3})$ (Theorem \ref{finthr}). We find $10$ generators: the classes $\gamma_{\Gamma}$ for all trees $\Gamma$ with at most three nodes (they are $5$), plus the Mumford class $\km_{2}$. The remaining $4$ generators are somewhat unexpected.

Unfortunately there are two problems in going beyond the case of three nodes, the one mentioned above with the Mumford classes, and the fact that the inductive technique will break down for curves with five nodes, because the top Chern class of the normal bundle to $\M_{0}^{5}$ in $\M_{0}^{\leq 5}$ is a $0$-divisor in $A^*(\M^5_0) \otimes \Q$ (Remark \ref{zerodiv}). 

{\bf Acknowledgments.} I am grateful to my thesis advisor Angelo Vistoli for his patient and constant guidance.

I also wish to thank Dan Edidin for helpful remarks and Rahul Pandharipande who suggested the problem to my advisor.

\newpage

\tableofcontents

\section{Description of the stack $\M_0$}

Let $\cat$ be the site of schemes over the complex point $\spc$ equipped with the \'etale topology.

\begin{defi}
We define the category of rational nodal curves $\M_0$ as the category over $\cat$ whose objects are flat and proper morphisms of finite presentation $C \xrightarrow{\pi} T$ (where $C$ is an algebraic space over $\C$ and $T$ is an object in $\cat$) such that for every geometric point $\spk \xrightarrow{t} T$ (where $\Omega$ is an algebraically closed field) the fiber $C_{t}$
\begin{equation*}
\xymatrix{
C_{t} \cart \ar[r]^{t'} \ar[d]_{\pi'} &C \ar[d]^{\pi}\\
\spk \ar[r]^t &T
}
\end{equation*}
is a projective reduced nodal curve such that
\begin{eqnarray*}
h^0(C_t, \Of_{C_t}) &=& 1 \; \text{(that is to say $C_t$ is connected;)}\\
h^1(C_t, \Of_{C_t}) &=& 0 \; \text{(thus the arithmetic genus of the curve is 0).}
\end{eqnarray*}
Morphisms in $\M_0$ are cartesian diagrams
\begin{equation*}
\xymatrix{
C' \cart \ar[r] \ar[d]_{\pi'} &C \ar[d]^{\pi}\\
T' \ar[r] &T
}
\end{equation*}
The projection $pr : \M_0 \xrightarrow{} \cat$ is the forgetful functor
\begin{eqnarray*}
pr:
&
 \left \{
\vcenter {\xymatrix{
C' \cart  \ar[r] \ar[d]_{\pi'} & C \ar[d]^{\pi}\\
T' \ar[r] & T
}}
\right \}
&
\mapsto \{ T' \xrightarrow{} T \}
\end{eqnarray*}
\end{defi}

We refer to \cite[Proposition 1.10]{fulg} for the proof that $\M_0$ is an Artin stack  in the sense of
\cite{Artd} (the proof uses standard arguments).

Here we describe more precisely the geometric points of $\M_0$. Given an algebraically closed field $\Omega$, let us define a {\it rational tree (on $\Omega$)} to be a connected nodal curve with finite components each of them is isomorphic to $\Pro^1_\Omega$ and which has no closed chains.
Classical cohomology argument (\cite[Proposition 1.3]{fulg}) allows us to state that geometric points of $\M_0$ are rational trees.

\begin{oss}It is important to notice that by definition geometric fibers of a family of curves $C \xrightarrow{\pi}T$ in $\M_0$ are projective and therefore are schemes. Moreover it is further known (see \cite {Kn} V Theorem 4.9) that  a curve (as an algebraic space) over an algebraically closed field $\Omega$ is a scheme.

Notwithstanding the fact that we consider families of curves whose fibers and bases are schemes we still need algebraic spaces because, as we show below, there exist families of rational nodal curves $C \xrightarrow{\pi} T$ where $C$ is not a scheme.
\end{oss}

\begin{exm} \label{HV}
Let us consider the following example which is based on Hironaka's example \cite{Hir} of an analytic threefold which is not a scheme (this example can be found also in \cite{Kn}, \cite{hrt} and \cite{Shf}).

Let $S$ be a projective surface over $\C$ and $i$ an involution without fixed points. Suppose that there is a smooth curve $C$ such that $C$ meets transversally the curve $i(C)$ in two points $P$ and $Q$. Let $M$ be the product $S \times \Pro ^1_{\C}$ and let $e: S \to M$ be the embedding $S \times 0$. On $M - e(Q)$, first blow up the curve $e(C - Q)$ and then blow up the strict transform of $e(i(C) - Q)$. On $M - e(P)$, first blow up the curve $e(i(C) - P)$ and then blow up the strict transform of $e(C - P)$. We can glue these two blown-up varieties along the inverse images of $M - e(P) - e(Q)$. The result is a nonsingular complete scheme $\widetilde M$. 

We have  an action of $\Cg_2$ (the order-2 group generated by the involution $i$) on $M$ induced by the action of $\Cg_2$ on $S$ and the trivial action on $\Pro^1_{\C}$. Furthermore we can lift the action of $\Cg_2$ to $\widetilde M$ such that we obtain an equivariant map
$$
\widetilde M \xrightarrow{\pi} S;
$$
we still call $i$ the induced involution on $\widetilde M$. We have an action of $\Cg$ on $\widetilde M$ which is faithful and we obtain geometric quotients $f$ and $g$
\begin{equation*}
\xymatrix{
\widetilde M \ar[d]_{\pi} \ar[r]^{f} & \widetilde M /\Cg_2  \ar[d]^{\widetilde{\pi}} \\
S \ar[r]^{g} & S/\Cg_2
}
\end{equation*}

Now  we notice that $\pi$ is a family of rational nodal curves. In particular the generic geometric fiber is isomorphic to $\Pro^1_{\C}$. Geometric fibers on the two curves $C$ and $i(C)$ have one node except those on $P$ and $Q$ which have two nodes. Similarly we have that also $\widetilde{\pi}$ is a family of rational nodal curves whose fibers with a node are on $g(C)=g(i(C))$ and there is only one fiber on $g(P)=g(Q)$.

Since $\Cg_2$ acts faithfully on $S$ and $S$ is projective we have that $S/\Cg_2$ is a scheme. Furthermore  we mention (see \cite{Kn} Chapter 4) that the category of separated algebraic spaces  is stable under finite group actions. Consequently the family $\widetilde{\pi}$ belongs to $\M_0(S/ \Cg_2)$. But $\widetilde M / \Cg_2$ is not a scheme (to prove this we can follow the same argument of \cite{hrt} Appendix B Example 3.4.2).

In order to complete our example, we have to exhibit a surface which satisfies the required properties. Take the Jacobian $J_2=\text{Div}^0(C)$ of a genus 2 curve $C$. 
Let us choose on $C$ two different points $p_0$ and $q_0$ such that $2(p_0 - q_0) \sim 0$. Let us consider the embedding
\begin{eqnarray*}
C &\to& J_2\\
p \in C &\mapsto& |p-q_0|.
\end{eqnarray*}
We still call $C$ the image of the embedding. Let $i$ be the translation on $J_2$ of $|p_0 - q_0|$. By definition $i$ is an involution such that $C$ and $i(C)$ meets transversally in two points: $0$ and $|p_0 - q_0|$.
\end{exm}

\subsection{Stratification of $\M_0$ by nodes} \label{strat}

For every object $C \xrightarrow{\pi} T$ in $\M_0$, we are going to define the {\it relative singular locus of $C$}, which will be denoted as $C_{rs}$: roughly speaking it is the subfamily of $C$ whose geometric fibers have nodes. Its image to the base $T$ is a closed subscheme. Moreover it is possible to define on $T$ closed subschemes $\{ T^{\geq k} \to T\}_{k \geq 0}$ where for each integer $k \geq 0$ the fiber product $T^{\geq k} \times_T C$ is a family of rational curves with at least $k$ nodes.
In order to give a structure of closed subspaces to $C_{rs}$ and to the various subschemes $T^{\geq k}$, we will describe them through Fitting ideals.

\begin{defi} We construct the closed subspace $C_{rs}$ and the $T^{\geq k}$ locally in the Zariski topology, so we consider a family of rational nodal curve $C \xrightarrow{\pi} \spe A$ for some $\C-$algebra $A$.

Now we follow \cite{mumlect} Lecture 8.

As the map $\pi$ has relative dimension 1 the relative differential sheaf $\Omega_{C/T}$ has rank 1 as a sheaf over $C$. Further because $\pi$ is a map of finite presentation and $A$ is a Noetherian ring, we have that the sheaf $\Omega_{C/T}$ is coherent. For every point $p \in T$ we set
 \begin{eqnarray*}
e(p)=  \dim _{k(p)} (\Omega _{C/T,p} \otimes k(p))
\end{eqnarray*}
where $k(p)$ is the residual field of $p$.

Choose a basis $\{a_i\}_{i=1, \dots, e(p)}$ of $\Omega _{C/T,p} \otimes k(p)$, they extend to a generating system for $\Omega _{C/T,p}$. Furthermore we have an extension of this generating system to $\Omega_{C/T}$ restricted to an \'etale  neighborhood of $p$. So we have a map
$$
\Of_{C} ^{\oplus e(p)} \to \Omega_{C/T}
$$
which is surjective up to a restriction to a possibly smaller neighborhood of 
$p$. At last (after a possibly further restriction) we have the following exact sequence of sheaves
\begin{eqnarray*}
\Of_{C}^{\oplus f} \xrightarrow{X} \Of_{C} ^{\oplus e(p)} \to \Omega_{C/T} \to 0
\end{eqnarray*}
where $X$ is a suitable matrix $f \times e$ of local sections of $\Of_{C}$. Let us indicate with $F _i(\Omega_{C/T}) \subset \Of_C$ the ideal sheaf generated by rank $e(p) - i$ minors of $X$. These sheaves are known as  {\it Fitting sheaves}
and they don't depend on the choice of generators (see \cite{lan} XIX, Lemma 2.3).

We define the {\it relative singular locus} to be the closed algebraic subspace $C_{rs} \subseteq C$ associated with $F_1(\Omega_{C/T})$.
\end{defi}

\begin{oss}Let us fix a point $t \in T$, we have that $t$ belongs to $T^{\geq k}$ iff
$$
F_{k-1}(\pi_*(\Of_{C_{rs}})_t) = 0 \Longleftrightarrow \dim_{k(t)} (\pi_*(\Of_{C_{rs}})_t \otimes k(t)) \geq k
$$
so $T^{\geq k}$ is the subscheme whose geometric fibers have at least $k$ nodes.
\end{oss}

\begin{defi}Set
\begin{eqnarray} \label{stratk}
T^0 &:=& T - T^{\geq 1}\\
T^k &:=& T^{\geq k} - T^{\geq k + 1}
\end{eqnarray}
\end{defi}

\begin{oss}We have that for every $k \in \N$ the subscheme $T^k$ is locally closed in $T$, further, given a point $t \in T$ and set 
$$
k:= \dim_{k(t)} (\pi_*(\Of_{C_{rs}})_t \otimes k(t)) \geq 0
$$
we have that $t$ belongs to a unique $T^k$.

Consequently, given a curve $C \xrightarrow{\pi} T$ where $T$ is an affine Noetherian scheme, the family 
$$
\{ T^k | T^k \neq \emptyset \}
$$
defined above is  {\it a stratification for $T$}.
\end{oss}

Now let $\St$ be an Artin stack over $\cat$. A family of locally closed substacks
$\{ \St^{\alpha}\}_{\alpha \in \N}$ represents a stratification for $\St$ if, for all morphisms
$$
T \to \St
$$
where $T$ is a scheme, the family of locally closed subschemes
$$
\{ T^{\alpha} := \St^{\alpha} \times_{\St} T \neq \emptyset \}
$$
is such that every point $t \in T$ is in exactly one subscheme $T^{\alpha}$, that is to say that $\{ T^{\alpha}\}$ is a stratification for $T$ in the usual sense.

\begin{defi}

We define $\M_0 ^{\geq k}$ as the full subcategory of $\M_0$ whose objects $C \xrightarrow{\pi} T$ are such that $T^{\geq k} = T$.

We further define $\M_0 ^{k}$ to be the subcategory whose objects $C \xrightarrow{\pi} T$ are such that $T^{k} = T$.

\end{defi}

\begin{oss}
We have from definition that for every morphism $T \to \M_0$
$$
T^{\geq k}=\M^{\geq k}_0 \times_{\M_0} T
$$
consequently $\{ \M^{k}_0 \}_{k \in \N}$ is a stratification for $\M_0$.
\end{oss}

\begin{propos} \label{stratreg}
For each ${k \in \N}$ the morphism
$$
\M_0 ^{\geq k} \xrightarrow{} \M_0
$$
is a regular embedding of codimension $k$.
\end{propos}

\begin{dem}
Let us consider a smooth covering
$$
U \xrightarrow{f} \M_0
$$
and let $C \xrightarrow{\pi} U$ be the associated curve.
we have to prove that the morphism
$$
U^k= \M_0 ^k \times_{\M_0} U \to U
$$
is a regular embedding of codimension $k$ (we already know that $U^k$ is a closed embedding).

Fix a point $p \in U^k\subset U$ and let $A:=\widehat{\Of}_{U,p} ^{sh}$ be the completion of the strict henselisation of the local ring $\Of_{U,p}$. (see \cite{egaq} Definiton 18.8.7). Consider the following diagram
\begin{eqnarray*}
\xymatrix{
C_A \cart \ar[d]_{\pi} \ar[r] & C \ar[d]_{\pi}\\
\spe A \ar[r] & U
}
\end{eqnarray*}
As $p$ is a point of $U^k$ we have that $C_{A, rs}$ is the union of $k$ nodes $q_1, \dots, q_k$. For each $i=1, \dots, k$ we have
$$
\widehat{\Of}_{C_A, q_i} =A \llbracket x,y \rrbracket /(xy-f_i)
$$
with $f_1, \dots , f_k \in A$.
So we can write
$$
M:=\Of_{C_{A,rs}} = \prod_{i=1}^k (A/f_i).
$$
Let us consider $M$ as an $A-$module, we have the following exact sequence of $A-$module
$$
A^k \xrightarrow{D} A^k \xrightarrow{} M \to 0
$$
where $D$ is the diagonal matrix
\begin{equation}
\left(\begin{array}{ccc}
f_1 & \cdots & 0 \\
\vdots & \ddots & \vdots \\
0 & \cdots & f_k
\end{array}\right)
\end{equation}
So we have $F_{k-1}(M)=(f_1, \dots, f_k)$ and this means that
$$
B:=\widehat{\Of}^{sh}_{U^{\geq k},p} = A/(f_1, \dots, f_k).
$$
From deformation theory we have that $\{ f_1, \dots, f_k \}$ is a regular sequence, consequently the map $U^{\geq k} \to U$ is regular of codimension $k$ as claimed.
\end{dem}

\subsection{Combinatorical version of rational nodal curves.}
Now let us fix a useful notation. Given an algebraically closed field $\Omega$ and an isomorphism class of $C$ (still denoted with $C$) in $\M_0(\spk)$, we define the {\it dual graph of $C$}, denoted $\Gamma(C)$ or simply $\Gamma$, to be the graph which has as many vertices as the irreducible components of $C$ and two vertices are joined by an edge if and only if the two corresponding lines meet each other.

For example we have the following correspondence
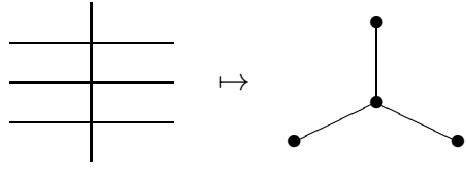
\begin{eqnarray*}
\vcenter{\trebg} & \; \mapsto \; & \vcenter{\treb}
\end{eqnarray*}
We can associate at least one curve with every tree by this map, but such a curve is not in general unique up to isomorphism; an example is the following tree
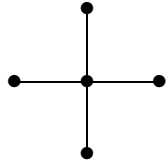
\begin{equation*}
\quattrob
\end{equation*}
In the following we think of a tree $\Gamma$ as a finite set of vertices with a connection law given by the set of pairs of vertices which corresponds to the edges.

\begin{defi}
Given a graph $\Gamma$ we call {\it multiplicity} of a vertex $P$ the number $e(P)$ of edges to which it belongs and we call $E(P)$ the set of  edges to which it belongs. Furthermore we call the {\it maximal multiplicity} of the graph $\Gamma$ the maximum of multiplicities of its vertices. We also call $\Delta_n$ the set of vertices with multiplicity $n$ and $\delta_n$ the cardinality of $\Delta_n$.
\end{defi}

\begin{oss} \label{isclass}
We can associate an unique isomorphism class of curves in $\M_0(\spk)$ to a given tree  $\Gamma$ if and only if the maximal multiplicity of $\Gamma$ is 3.
\end{oss}

Tree graphs classify the topological type of rational nodal curves.

Now let us fix a stratum $\M_0^k$, there are as many topological types of curves with $k$ nodes as trees with  $k+1$ vertices.

Purely topological arguments show the following

\begin{lemma} Given a curve $C \xrightarrow{\pi} T$ in $\M_0^{k}$ where $T$ is a connected scheme, the curves of the fibers are of the same topological type. 
\end{lemma}

So we can give the following:

\begin{defi}
For each tree $\Gamma$ with $k+1$ vertices we define $\M_0^{\Gamma}$ as the full subcategory of $\M_0 ^k$ whose objects are curves $C \xrightarrow{\pi} T$ such that on each connected component of $T$ we have curves of topological type $\Gamma$.
\end{defi}

For each $\Gamma$, $\M^{\Gamma}_0$ is an open (and closed) substack of $\M^{k}_0$ and we can consequently write
$$
\M_0 ^k = \coprod_{\Gamma} \M_0 ^{\Gamma}
$$
where $\Gamma$ varies among trees with $k+1$ vertices.

\subsection{Description of particular strata}\label{strata}

From now on we focus on graphs (and curves) with maximal multiplicity equal to 3. 

\begin{lemma}
Let $\Gamma$ be a graph with maximal multiplicity equal to 3 and $C$ an isomorphic class of curves of topological type $\Gamma$, then we have the equivalence
\begin{equation} \label{equiv}
\M_0 ^{\Gamma} \simeq \B \aut(C).
\end{equation}
\end{lemma}

\begin{dem}
From Remark (\ref{isclass}) we have that all curves of the same topological type $\Gamma$ are isomorphic.
\end{dem}

We have a canonical surjective morphism
$$
\text{Aut}(C) \xrightarrow{g} \text{Aut}(\Gamma)
$$
which sends each automorphism to the induced graph automorphism.

\begin{propos} \label{coord}
There is a (not canonical) section $s$ of $g$
\begin{equation*}
\xymatrix{
\text{Aut}(C) \ar@<-1ex>[r]_g &  \text{Aut}(\Gamma) \ar@<-1ex>[l]_s 
}
\end{equation*}
\end{propos}

\begin{dem}
Let us fix coordinates $[X,Y]$ on each component of $C$ such that
\begin{itemize}
\item on components with one node the point $[1,0]$ is the node,
\item on components with two nodes the points $[0,1]$ and $[1,0]$ are the nodes,
\item on components with three nodes the points $[0,1]$, $[1,1]$ and $[1,0]$ are the nodes.
\end{itemize}
Let us define on each component
\begin{eqnarray*}
0:=[0,1] \quad 1:=[1,1] \quad \infty:=[1,0]
\end{eqnarray*}

In order to describe the section $s$, let us notice that, given an element $h$ of $\aut(\Gamma)$, there exists a unique automorphism $\gamma$ of $C$ such that:
\begin{itemize}
\item $\gamma$ permutes components of $C$ by following the permutation of vertices given by $h$
\item on components with one node, $\gamma$ makes correspond the points $0,1$
\item on components with two nodes, $\gamma$ makes correspond the point $1$.
\end{itemize}
\end{dem}

\begin{propos}
We have
\begin{equation*}
\aut (C) \cong \aut(\Gamma) \ltimes \left( (\Gm)^{\Delta_2} \times E^{\Delta_1} \right).
\end{equation*}
where $E$ is the subgroup of $\Pro GL_2$ that fixes $\infty$ and $\Gm$ the multiplicative group of the base field.
\end{propos}

\begin{dem}
Let us consider the normal subgroup $g^{-1}(\id)$ of automorphisms of $C$ which do not permute components. It is the direct product of groups of automorphisms of each component that fixes nodes.

On components with one node the group of automorphisms is the subgroup $E$ of $\Pro GL_2$ that fixes $\infty$. $E$ can be described as the semidirect product of $\Gm$ and $\Ga$ (the additive group of the base field), moreover, having fixed coordinates, we can have an explicit split sequence
$$
\xymatrix{
0 \ar[r] & \Ga  \ar[r]^{\varphi} & E \ar@<-1ex>[r]_{\rho} & \Gm \ar@<-1ex>[l]_{\psi} \ar[r] & 1.
}
$$

On components of two nodes we have that the group of automorphisms is $\Gm$.

At last only identity fixes three points in $\Pro^1_\Omega$. So we can conclude that
\begin{equation*}
g^{-1}(\id)=E^{\Delta_1} \times \Gm^{\Delta_2}
\end{equation*}
and we have a (not canonical) injection
\begin{equation*}
E^{\Delta_1} \times \Gm^{\Delta_2} \to \text{Aut}(C)
\end{equation*}

Then $Aut(C)$ is the semi-direct product given by the exact sequence
\begin{equation} \label{seps}
\xymatrix{
1 \ar[r] & \Gm^{\Delta_2} \times E^{\Delta_1} \ar[r] & \text{Aut}(C) \ar@<-1ex>[r] &  Aut(\Gamma) \ar@<-1ex>[l] \ar[r] & 1
}
\end{equation}
and we write it as
\begin{equation*}
\text{Aut}(\Gamma) \ltimes (\Gm)^{\Delta_2} \times E^{\Delta_1}.
\end{equation*}
\end{dem}

So we can explicit
\begin{equation} \label{ps}
\M^{\Gamma}_0= \B \left( \text{Aut}(\Gamma) \ltimes (\Gm)^{\Delta_2} \times E^{\Delta_1} \right).
\end{equation}

\subsection{Dualizing and normal bundles} \label{normal}

In this section we briefly give the description of basic bundles over $\M_0$ and its strata $\M^k_0$. In particular we point out their restriction to $\M_0^{\leq 3}$.

{\bf The dualizing bundle.} On $\M_0$ we consider the universal curve $\cu \xrightarrow{\Pi} \M^{\Gamma}_0$.
defined in the following way:
\begin{defi}
Let $\cu$ be the fibered category on $\cat$ whose objects are families $C \xrightarrow{\pi} T$ of $\M_0$ equipped with a section $T \xrightarrow{s}C$ and whose arrows are arrows in $\M_0$ which commute with sections.
\end{defi}

Let $U \to \M_0$ be a smooth covering ($\M_0$ is an Artin stack) and let us consider the following cartesian diagram

\begin{equation*}
\xymatrix{
 \cu_U \cart \ar[d]^\Pi \ar[r] & \cu \ar[d]^\Pi\\
U \ar[r] & \M_0
}
\end{equation*}

\begin{defi}
We define the relative dualizing sheaf $\omega_0$ of $\cu \xrightarrow{\Pi} \M_0$ as the dualizing sheaf $\omega_{\cu_U/ U}$ on $\cu_U$.
\end{defi}

This is a good definition because for each curve $C \xrightarrow{\pi} T$ in $\M_0$ there is  the dualizing sheaf $\omega_{C/T}$ and its formation commutes with the base change.

For each curve $C \xrightarrow{\pi} T$ in $\M_0^{\leq k}$ we have the push forward $\pi_* \omega^{\vee}_{C/T}$. If we wish to have a well defined push forward of $\left ( \omega_{0} \right )^{\vee}$ we need to prove that for each curve in $\M_0^{\leq k}$ the sheaf $\pi_{*} \omega ^{\vee}_{C/T}$ is locally free of constant rank and it respects the base change. We can do it when $k$ is 3. In particular we have the following

\begin{propos} \label{dualsh}
Let $C \xrightarrow{\pi} T$ be a curve in $\M_0^{\leq 3}$. Then $\pi_{*} \omega_{C/T}^{\vee}$ is a locally free sheaf of rank 3 and its formation commutes with base change.
\end{propos}

\begin{dem}
First of all we prove that it is locally free and its formation commutes with base change.
It is enough to show that for every geometric point $t$ of $T$ 
$$
H^1(C_t, \omega^{\vee}_{C_t/t})=0.
$$
From Serre duality
$$
H^1(C_t, \omega^{\vee}_{C_t/t})=H^0(C_t, \omega^{\otimes 2}_{C_t/t})^{\vee}.
$$
When the fiber is isomorphic to $\Pro ^1$ we have 
$$
\omega^{\otimes 2}_{C_t/t} = (\Of(-2))^{\otimes 2} = \Of(-4)
$$
and we do not have global sections different from 0.

When $C_t$ is singular we have that the restriction of $\omega^{\otimes 2}_{C_t/t}$ to components with a node is $(\Of(-1))^{\otimes 2}=\Of(-2)$ and consequently the restriction of sections to these components must be 0. The restriction to components with two nodes is $\Of(0)$ (so restriction of sections must be constant) and restriction to components with three nodes is $\Of(2)$ (in this case the restriction of sections is a quadratic form on $\Pro^1$).

Given these conditions, global sections of $\omega^{\otimes 2}_{C_t/t}$ on curves with at most three nodes (as they have to agree on the nodes)
must be zero.

Similarly we verify that $h^0(C_t, \omega^{\vee}_{C_t/t}) = 3$ and conclude by noting that
$$
(\pi_*(\omega^{\vee}_{C/T}))_t = H^0 (C_t, \omega^{\vee}_{C_t/t}).
$$
\end{dem}

We notice that in $\M_0^{\leq 4}$ there are curves $C \xrightarrow{\pi} \spk$ for which there are global sections for $\omega^{\otimes 2}_{C/\spk}$. When we have a curve of topological type
$$
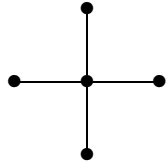
\quattrob
$$
the restriction of global sections on the central component are quartic forms on $\Pro^1_{\Omega}$ which vanishes on four points, thus we have
$$
H^0(C, \omega^{\otimes 2}_{C/\spk}) = \Omega.
$$ 

In conclusion we have a well defined bundle $\Pi_*\left ( \omega_{0} \right )^{\vee}$ on $\M^{\leq 3}_0$ and we refer to it as the dualizing bundle of $\M^{\leq 3}_0$.

\medskip

{\bf The normal bundles.} Given a tree $\Gamma$ we consider the (local) regular embedding (see Proposition \ref{stratreg})
$$
\M_0^\Gamma \xrightarrow{in} \M_0^{\leq \delta}
$$
where $\delta$ is the number of edges in $\Gamma$.
From \cite{kre} Section 5 we have that there exists a relative tangent bundle $\mathcal T_{in}$ on $\M_0^\Gamma$ which injects in $in^*(\mathcal T_{\M_0^{\leq \delta}})$.
\begin{defi}
Let us consider the following exact sequence of sheaves
$$
0 \to \mathcal T_{in} \to in^*(\mathcal T_{\M_0^{\leq \delta}}) \to in^*(\mathcal T_{\M_0^{\leq \delta}}) / \mathcal T_{in} \to 0.
$$
We define the normal bundle as the quotient sheaf on $\M_0^\Gamma$
$$
\Nf_{in}:=in^*(\mathcal T_{\M_0^{\leq \delta}}) / \mathcal T_{in}
$$
\end{defi}
When $\Gamma$ has maximal multiplicity at most 3, we can describe it as the quotient of the space of first order deformations by $\aut(C)$ in the following way.
Let us consider the irreducible components $ C_{\Gamma}$ of $C$. The space of first order deformations near a node $P$ where two curves $C_\alpha$ and $C_\beta$ meet is (see \cite{HM} p. 100)
$$
T_P(C_\alpha) \otimes T_P(C_\beta).
$$
Consequently we have the following:

\begin{lemma}
The space $N_\Gamma$ of first order deformations of $C$ is
$$
{\bigoplus_{P \in E(\Gamma)}} T_P(C_\alpha) \otimes T_P(C_\beta).
$$
\end{lemma}

On $N_\Gamma$ there is an action of $\aut(C)$ which we will describe in Section \ref{stratclass}.

\subsection{Algebraic lemmas} Here we state and prove two algebraic lemmas for future reference. 

The first one is Lemma 4.4 \cite{VV}:
\begin{lemma} \label{incollamento}
Let $A, B$ and $C$ be rings, $f:B \to A$ and $g:B \to C$ ring homomorphism. Let us suppose that there exists an homomorphism of abelian groups $\phi: A \to B$ such that the sequence
$$
A \xrightarrow{\phi} B \xrightarrow{g} C \xrightarrow{} 0
$$
is exact; the composition $f \circ \phi: A \to A$ is multiplication by a central element $a \in A$ which is not a 0-divisor.

Then $f$ and $g$ induce an isomorphism of ring
$$
(f,g):B \to A \times_{A/(a)} C,
$$
where the homomorphism $A \xrightarrow{p} A/(a)$ is the projection, while $C \xrightarrow{q} A/(a)$ is induced by the isomorphism $C \simeq B/\ker g$ and the homomorphism of rings $f: B \to A$.
\end{lemma}

\begin{dem}
Owing to the fact that $a$ is not a 0-divisor we immediately have that $\phi$ is injective. Let us observe that the map $(f,g):B \to A \times_{A/(a)} C$ is well defined for universal property and for commutation of the diagram:
\begin{equation*}
\xymatrix{
B  \ar[r]^g \ar[d]_{f} &C \ar[d]^{q}\\
A \ar[r]^{p} &A/(a)
}
\end{equation*}
Now let us exhibit the inverse function $A \times_{A/(a)} C \xrightarrow{\rho} B$. Given $(\alpha, \gamma) \in  A \times_{A
/(a)} C$, let us chose an element $\beta \in B$ such that $g(\beta)=\gamma$. By definition of $q$ we have that $f(\beta)-\alpha$ lives in the ideal $(a)$ and so (it is an hypothesis  on $f \circ \phi$) there exists in $A$ an element $\tilde \alpha$ such that:
$$
f(\beta)-\alpha=f(\phi(\tilde \alpha))
$$
from which:
$$
f(\beta - \phi( \tilde \alpha)) = \alpha
$$
we define then $\rho (\alpha,\gamma):=\beta-\phi( \tilde \alpha)$. In order to verify that it is a good definition, let us suppose that there exist an element $\beta_0 \in B$  such that  $(f,g)(\beta_0)=0$, to be precise there exists an element $\alpha_0 \in A$ such that: $\phi(\alpha_0)=0$ and furthermore
$$
0=f(\phi(\alpha_0))=a\alpha_0
$$
but we have that $a$ is not a divisor by zero, so necessary we have $\alpha_0=0$ and $\beta_0=0$. we conclude by noting that from the definition of $\rho$ we have immediately that it is an isomorphism. 
\end{dem}

\begin{oss}
The Lemma \ref{incollamento} will be used for computing the Chow ring of $\M_0^{\leq k+1}$ when the rings $A^*[\M_0^{\leq k}]\rat$ and $A^*[\M_0^{k+1}] \rat$ are known. As a matter of fact, given an Artin stack $\mathcal X$ and a closed Artin substack $\mathcal Y \xrightarrow{i} \mathcal X$ of positive codimension, we have the exact sequence of groups (see \cite{kre} Section 4)
$$
A^*(\mathcal Y) \rat \xrightarrow{i_*} A^*(\mathcal X) \rat \xrightarrow{j^*} A^*(\mathcal U) \xrightarrow{} 0
$$
By using the Self-intersection Formula, it follows that: 

$$
i^*i_*1=i^*[\mathcal Y]=c_{top}(\Nf_{\mathcal Y / \mathcal X}),
$$
when $c_{top}(\Nf_{\mathcal Y / \mathcal X})$ is not 0-divisor we can apply the Lemma.
\end{oss}

We will also use the following algebraic Lemma:

\begin{lemma}\label{quoz}
Given the morphisms

$$
A_1 \xrightarrow{p_1} \overline{A}_1 \xrightarrow{} B \xleftarrow{} \overline{A}_2 \xleftarrow{p_2} A_2
$$
in the category of rings, where the maps $p_1$ and $p_2$ are quotient respectively for ideals  $I_1$ and $I_2$. Then it defines an isomorphism
$$
\overline{A}_1 \times_B \overline{A}_2 \cong \frac{A_1 \times_B A_2}{(I_1,I_2)}.
$$
\end{lemma}
\begin{dem}
Let us consider the map
$$
(p_1, p_2): A_1 \times_B A_2 \to \overline{A}_1 \times_B \overline{A}_2,
$$
by surjectivity of $p_1$ e $p_2$ this map is surjective, while the kernel is the ideal $(I_1,I_2)$.
\end{dem}

\section{Fundamental classes on $\M_0$}

In this Section we work on a fixed tree $\Gamma$ with maximal multiplicity k and $\delta$ edges. We indicate with $\Delta_1, \Delta_2, \dots, \Delta_k$ the sets of vertices that belongs respectively to one, two and three edges (and with $\delta_1, \delta_2, \dots, \delta_k$ their cardinalities).

\subsection{Chow rings of strata}\label{stratclass}

We restrict the universal curve $\cu \to \M_0$ to $\M^{\Gamma}_0$
$$
\cu^{\Gamma} \xrightarrow{\Pi} \M^{\Gamma}_0.
$$
and consider the normalization (see \cite{Vis} Definition 1.18)
$$
\widehat{C}^{\Gamma} \xrightarrow{N} \cu.
$$

\begin{oss}
Given a curve $C \xrightarrow{\pi} T$ in $\M^{\Gamma}_0$ (that is to say a morphism $T \xrightarrow{f} \M^{\Gamma}_0$) we define $\widehat{C}$ as $T \times_{\M^{\Gamma}_0} \widehat{\cu}^{\Gamma}_0$. We have the following cartesian diagram
\begin{equation*}
\xymatrix{
\widehat{C} \cart \ar[r] \ar[d]_{n} & \widehat{\cu}^{\Gamma} \ar[d]^{N}\\
C \cart \ar[d]_{\pi} \ar[r] & \cu \ar[d]^{\Pi}\\
T \ar[r]^{f} & \M^{\Gamma}_0.
}
\end{equation*}

We notice that when $T$ is a reduced and irreducible scheme $\widehat{C} \xrightarrow{n} C$ is the normalization.

The map $\pi n: \widehat{C} \to T$ is proper as $\Pi N: \cu^{\Gamma} \to \M^{\Gamma}_0$ is. Then there exists a finite covering $\widetilde T \xrightarrow{} T$ (see \cite{Kn} Chapter 5, Theorem 4.1) such that we have the following commutative diagram
\begin{equation*}
\xymatrix{
 & \widehat C \ar[ddl]_g \ar[d]^n\\
 & C \ar[d]^\pi\\
\widetilde{T} \ar[r] &T
}
\end{equation*}
where the map $g$ has connected fibers.
\end{oss}

\begin{defi}
We define $\widetilde{\M}_0 ^{\Gamma}$ as the fibered category on $\cat$ whose objects are rational nodal curves $C \xrightarrow{\pi} T$ in $\M_0^{\Gamma}$ equipped with an isomorphism $\varphi: \coprod^{\Gamma}T  \to \widetilde T$ over $T$ and maps are morphisms in $\M_0^{\Gamma}$ that preserve isomorphisms.
\end{defi}

Straightforward arguments show the following

\begin{lemma}
The forgetful morphism
$$
\widetilde{\M}_0 ^{\Gamma} \to \M_0^{\Gamma}
$$
is representable finite \'etale and surjective.
\end{lemma}

We have a more explicit description of $\widetilde{\M}_0 ^{\Gamma}$: if we call $\Ms_{0,i}^n$ the stack of rational curves with $n$ nodes and $i$ sections, we can exhibit an equivalence
\begin{equation} \label{equivstrat}
\widetilde{\M}_0^{\Gamma} \cong \prod_{\alpha \in \Gamma} (\Ms_{0,e(\alpha)}^0)
\end{equation}
where, for each $\alpha \in \Gamma$, $e(\alpha)$. The proof of this is straightforward, anyway in the following we assume that the maximal multiplicity $k$ of $\Gamma$ is at most 3. In this case, as we have seen in  Section \ref{strata}:
$$
\M_0^{\Gamma} = \B \aut(C)
$$
where $\aut(C)$ is the group $\aut(\Gamma) \ltimes \left ( E^{\Delta_1} \times \Gm^{\Delta_2} \right )$. Let us call $\sigma$ the order of $\aut{\Gamma}$.

The \'etale covering of degree $\sigma$
\begin{eqnarray*}
\widetilde{\M}_0^{\Gamma} &\xrightarrow{\phi}& \M_0^{\Gamma}
\end{eqnarray*}
becomes
\begin{eqnarray*}
\B H &\xrightarrow{\phi}& \B \aut(C)
\end{eqnarray*}
where $H$ is the group $E^{\Delta_1} \times \Gm^{\Delta_2}$.

In order to compute the Chow ring on each stratum we only need equivariant intersection theory. 

\begin{propos} \label{ringstratum}
Let $\Gamma$ be a tree of maximal multiplicity at most 3 and let $C$ be the unique isomorphism class of curves of topological type $\Gamma$. Let us fix on $C$ coordinates as in the proof of Proposition \ref{coord}.

Then the Chow ring $A^*(\M^{\Gamma}_0)$ is
$$
A^*_{\aut(C)} \rat \cong (\Q[t_{\Delta_1}, r_{\Delta_2}])^{\aut(\Gamma)}
$$
where the action of an element $g \in \aut(\Gamma)$ on $\Q[x_{\Delta_1}, y_{\Delta_2}]$ is the obvious permutation on the $\Delta_1 \cup \Delta_2$ variables together with multiplication by (-1) of $r$-variables corresponding to components of which $g$ exchanges $0$ and $\infty$.
\end{propos}

\begin{dem}
The group $E^{\Delta_1} \times \Gm^{\Delta_2}$ is a normal subgroup of $\aut(C)$.

So we can apply the following

\begin{lemma} \label{isophi} \cite{Ve}
Given an exact sequence of algebraic groups over $\C$
$$
\xymatrix{
1 \ar[r] & H  \ar[r]^\phi & G \ar[r]^\psi &  F \ar[r] & 1.
}
$$
with $F$ finite and $H$ normal in $G$, we have
$$
A^*_{G} \rat \cong (A^*_{H} \rat )^{F}
$$
\end{lemma}
and conclude that
$$
A^*_{\text{Aut}(C)} \rat \cong (A^*_{ E^{\Delta_1}\times (\Gm)^{\Delta_2}} \rat)^{\text{Aut}(\Gamma)}
$$

So we have reduced to compute $A^*_{E^{\Delta_1} \times (\Gm)^{\Delta_2}}$ and the action of $\aut(\Gamma)$ on it.

\begin{claim} \label{pr}
The ring $A^*_{E^{\Delta_1} \times (\Gm)^{\Delta_2} } \rat$ is $\mathbb Q[x_{\Delta_1}, y_{\Delta_2}]$, that is to say it is algebraically generated by  $\Delta_1 \cup \Delta_2$ independent generators of degree 1.
\end{claim}

We use the following fact (\cite{Ve} Proposition 2.8):
\begin{lemma}
For every linear algebraic group $G$, we have
$$
A^*_{G \times \Gm} \cong A^*_{G} \otimes_{\Z} A^*_{\Gm}.
$$
\end{lemma}
By recalling that $ A^*_{\Gm} \rat$ is isomorphic to $\Q [r]$, where $r$ is an order one class, we have
$$
A^*_{E^{\Delta_1} \times (\Gm)^{\Delta_2} } \rat \cong A^*_{E^{\Delta_1}} \otimes_{\Q} \Q[y_{\Delta_2}].
$$
Using the fact that the group $E$ is the semidirect product
$$
\xymatrix{
0 \ar[r] & \Ga  \ar[r]^{\varphi} & E \ar@<-1ex>[r]_{\rho} & \Gm \ar@<-1ex>[l]_{\psi} \ar[r] & 1,
}
$$
we can explicit that
$$
A^*_{E^{\Delta_1}} \rat \cong \Q[x_{\Delta_1}].
$$

Now we describe the action of $\aut(\Gamma)$ on
$$
A^*_{E^{\Delta_1} \times \Gm^{\Delta_2}} \cong \Q[x_{\Delta_1}, y_{\Delta_2}].
$$
In order to do this we need a more explicit description of the classes $x_{\Delta_1}, y_{\Delta_2}$ through the equivalence (\ref{equivstrat})

\begin{eqnarray*}
(\Ms^0_{0,1})^{\Delta_1} \times (\Ms^0_{0,2})^{\Delta_2} \times (\Ms^0_{0,3})^{\Delta_3} &\cong& \widetilde{\M}_0^{\Gamma}.
\end{eqnarray*}
Since $\Ms^0_{0,3} \simeq \spe \C$ we have
\begin{eqnarray*}
(\Ms^0_{0,1})^{\Delta_1} \times (\Ms^0_{0,2})^{\Delta_2} &\xrightarrow{\phi}& \M_0^{\Gamma}.
\end{eqnarray*}
Moreover we have
\begin{eqnarray*}
\Ms^0_{0,1} &\simeq& \B E\\
\Ms^0_{0,2} &\simeq& \B \Gm.
\end{eqnarray*}

Let $\alpha$ be a vertex of $\Gamma$ such that $e(\alpha)=1$ or $2$. On the component $\Ms_{0,e(\alpha)}^0$ let us consider the universal curve
$$
\cu^{\alpha} \xrightarrow{\widetilde{\Pi}} \Ms_{0,e(\alpha)}^0. 
$$
On $\Pro^1_{\C}$ we fix coordinates and we define the points
\begin{eqnarray*}
z_\infty &:=& [1,0]\\\
z_0 &:=& [0,1]
\end{eqnarray*}
We can write
\begin{itemize}
\item $\cu^{\alpha} \simeq [\Pro^1_{\C}/E]$ when $e(\alpha)=1$
\item $\cu^{\alpha} \simeq [\Pro^1_{\C}/\Gm]$ when $e(\alpha)=2$
\end{itemize}
Let us consider on $\Pro^1_{\C}$ the linear bundles $\Of(z_\infty)$ and $\Of(z_0)$. We have a natural action of $\Gm$ on global sections of both of them induced by the action of $\Gm$ (that, we recall, fixes $z_\infty$ and $z_0$) on $\Pro^1_{\C}$. Similarly we have an action of the group $E$ on global sections of $\Of(z_\infty)$.
When $e(\alpha)=1$ set
\begin{equation*}
\psi^{1}_{\infty,\alpha}:=c^E_1(H^0(\Of(z_\infty), \Pro^1_{\C}));
\end{equation*}
while, if $e(\alpha)=2$ set
\begin{eqnarray*}
\psi^{2}_{\infty,\alpha} &:=& c^{\Gm}_1(H^0(\Of(z_\infty), \Pro^1_{\C}))\\
\psi^{2}_{0,\alpha} &:=& c^{\Gm}_1(H^0(\Of(z_0), \Pro^1_{\C})).
\end{eqnarray*}
Clearly, when $e(\alpha)=2$, we have
$$
\psi^2_{\infty,\alpha} + \psi^2_{0,\alpha}=0
$$
Now we define
\begin{eqnarray*}
&& t_{\alpha}:=\psi^1_{\infty,\alpha} \text{ when $e(\alpha)=1$}\\
&& r_{\alpha}:=\psi^2_{\infty,\alpha} \text{ when $e(\alpha)=2$}
\end{eqnarray*}
We can write for each $\alpha$ such that $e(\alpha)=2$
$$
r_\alpha=\frac{\psi^2_{\infty,\alpha} - \psi^2_{0,\alpha} }{2}.
$$
Clearly all the classes $t_{\Delta_1}$ and $r_{\Delta_2}$ are of order one and independent. From what we have seen above these classes generates the ring $A^*(\widetilde{\M}^{\Gamma}_0) \rat$ and we have
$$
A^*(\widetilde{\M}^{\Gamma}_0) \rat \simeq \Q[t_{\Delta_1}, r_{\Delta_2}].
$$
Now we can describe the action of $\aut{\Gamma}$ on $\Q[t_{\Delta_1}, r_{\Delta_2}]$. An element $g \in \aut(\Gamma)$ acts on $C$ with a permutation $g_1$ on the components with one node and a permutation $g_2$ of components with two nodes. As we have chosen coordinates on $C_0$ such that $\infty$ corresponds to the node of the terminal components, we make $g_1$ act directly to the set $\{ t_{\Delta_1} \}$. We make $g_2$ act similarly on the set $\{ r_{\Delta_2} \}$ but we have in addition to consider the sign, that is to say that when $g_2$ sends a vertex $P$ of $\Gamma$ to another vertex $\beta$ (such that $e(\alpha)=e(\beta)=2$), we have two possibilities
\begin{itemize}
\item the automorphism $g$ exchange coordinates 0 and $\infty$ and so we have
$$
g(r_{\alpha})=g \left( \frac{\psi^2_{\infty,\alpha} - \psi^2_{0,\alpha}}{2}\right ) = \frac{\psi^2_{0,\beta} - \psi^2_{\infty,\beta}}{2} = -r_{\beta}
$$
\item the automorphism $g$ sends 0 in 0 and $\infty$ in $\infty$; in this case we have
$$
g(r_{\alpha}) = r_\beta.
$$
 \end{itemize}

\end{dem}

\begin{defi}
We define  $\gamma_i$ as the class in  $A^*(\M_0^{\leq i})$ of $\M_0^i$. We will indicate with $\gamma_i \in A^*(\M_0)$ also the class of the closure of $\M_0^i$ in $\M_0$. Similarly we define $\gamma_{\Gamma}$ as the class of the closure of $\M^{\Gamma}_0$ in $\M_0$.
\end{defi}

\begin{propos}\label{classstratum}
Let $\Gamma$ be a tree of maximal multiplicity at most three (except the single point) and $C$ be the curve of topological type $\Gamma$. Let us consider the \'etale covering
$$
\widetilde{\M}_0^\Gamma \xrightarrow{\phi} \M_0^{\Gamma}
$$
Let $C_{\Gamma}$ be the components of $C$ (which we see as vertex of $\Gamma$). Let $E(\Gamma)$ be the set of edges. If $(\alpha, \beta) \in E(\Gamma)$ we call $z_{\alpha \beta}$ the common point of $C_\alpha$ and $C_\beta$.
Then, by following notation of Proposition \ref{ringstratum}, we have
$$
\phi^* \gamma_{\Gamma} = \prod_{(\alpha, \beta) \in E(\Gamma)} \left( \psi^{e(\alpha)}_{z_{\alpha \beta}|_{\alpha}, \alpha} + \psi^{e(\alpha)}_{z_{\alpha \beta}|_{\beta}, \beta} \right).
$$
In particular the classes of each stratum of $\M^{\leq 3}_0$
after fixing coordinates on $C$ and ordering components of $\Delta_1$ and $\Delta_2$, are:
\medskip

\begin{center}
\begin{tabular}{|c|c|}
\hline
\small{Graph $(\Gamma)$} & class of stratum $\phi^*\gamma_{\Gamma}$\\
\hline
\small{$$ \unoa $$} & $t_1 + t_2$\\
\hline
\small{$$ \duea $$} & $(t_1 - r_1)(t_2 + r_1)$ \\
\hline
\small{$$ \trea $$} & $(t_1 - r_1)(r_1 + r_2)(t_2 - r_2)$\\
\hline
\small{$$ \treb $$} & $t_1 t_2 t_3$\\
\hline
\end{tabular}
\end{center}

\end{propos}

\begin{dem}
For every tree $\Gamma$ of maximal multiplicity at most three (except the single point), let us consider the regular embedding
$$
\M_0^\Gamma= \B \aut(C) \xrightarrow{in} \M_0^{\leq \delta}
$$
where $\delta$ is the number of edges of $\Gamma$.

Let us consider the normal bundle
$$
N_\Gamma:= N_{\M_0^\Gamma / \M_0^{\leq \delta}}.
$$
As we have recalled in Section \ref{normal} we have that $N_\Gamma=\text{def}_{\Gamma}$ is the space of first order deformations of $\M_0^{\Gamma}$
$$
{\bigoplus_{(\alpha,  \beta) \in E(\Gamma)}} T_{z_{\alpha \beta}}(C_\alpha) \otimes T_{z_{\alpha \beta}}(C_\beta).
$$

As usual let us consider the \'etale covering
\begin{eqnarray*}
&& \widetilde \M_0^{\Gamma} \xrightarrow{\phi} \M_0^{\Gamma}\\
&& \B H \to \B ( \aut (\Gamma) \ltimes H )
\end{eqnarray*}
and set $\widetilde N _{\Gamma}= \phi^* N_{\Gamma}$. We notice that, given coordinates as in Section \ref{strat}, the point $z_{\alpha\beta}$ on each component $C_\alpha$ (which we call $z_{\alpha\beta}|_{\alpha}$) is 0, 1 or $\infty$.

By using notation of Proposition \ref{ringstratum} we have on $\widetilde \M^{\Gamma}_0$
$$
c^{G_{\alpha}}_1(T_{z_{\alpha \beta}}(C_\alpha))= \psi^{e(\alpha)}_{z_{\alpha \beta}|_{\alpha}, \alpha}
$$
where $G_{\alpha}=E$ if $e(\alpha)=1$, $G_{\alpha}=\Gm$ if $e(\alpha)=2$ and $G_{\alpha}=\id$ if $e(\alpha)=3$
Notice that $c^{G_{\alpha}}_1(T_{z_{\alpha \beta}}$ is zero when $e(\alpha)=3$, consequently
$$
c^{H}_{\text{top}}(\widetilde{N}_{\Gamma})= \prod_{(\alpha, \beta) \in E(\Gamma)} \left( \psi^{e(\alpha)}_{z_{\alpha \beta}|_{\alpha}, \alpha} + \psi^{e(\beta)}_{z_{\alpha \beta}|_{\beta}, \beta} \right).
$$
We have the following relation in $A^*(\M^{\Gamma}_0)$

\begin{equation} \label{nor}
in^*in_*[\M_0^{\Gamma}]=c^{H}_{top}(N_\Gamma) \cap [\M_0^{\Gamma}].
\end{equation}
and so
$$
\phi^*\gamma_{\Gamma}= c^{H}_{\text{top}}(\widetilde{N}_{\Gamma})
$$
\end{dem}

\begin{oss} \label{zerodiv}

Given a choice of coordinates the class $\phi^* \gamma_{\Gamma}$ is invariant for the action of $\aut (\Gamma)$ given in Proposition \ref{ringstratum}, so we actually can see it as the class $\gamma_{\Gamma}$ in $\M^\Gamma_0$.

Moreover we have shown that these classes are not 0-divisor in $A^*(\M^{\Gamma}_0)$ for each $\Gamma$ corresponding to a stratum in $\M^{\leq 3}_0$, so we can apply Lemma \ref{incollamento}.

This fails if we consider more than four nodes. For example if we consider the following graph $\Gamma$
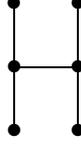
\begin{equation*}
\cinquenz
\end{equation*}
we have $\phi^* \gamma_{\Gamma}=0$.
\end{oss}

For future reference we can now state the following

\begin{propos} \label{injleqtre}
The Chow ring $A^*(\M^{\leq 3}_0) \rat$ injects into the product of $A^*(\M^{\Gamma}_0)\rat$ over trees with at most three edges.
\end{propos}
\begin{dem}
From Proposition \ref{classstratum} and Remark (\ref{zerodiv}), we have, for each $\delta \leq 3$, the following exact sequence of additive groups
$$
0 \to A^*(\M^{\delta}_0) \rat \xrightarrow{i^{\delta}_*} A^*(\M^{\leq \delta}_0) \rat \xrightarrow{j^{*\delta}} A^*(\M^{\leq (\delta -1)}_0)
$$
where
\begin{eqnarray*}
&&i^{\delta}: \M^{\delta}_0 \to \M^{\leq \delta}_0\\
&&j^{\delta}: \M^{\leq (\delta -1)}_0 \to \M^{\leq \delta}_0
\end{eqnarray*}
are the natural closed embeddings.

Let us consider the morphism
$$
A^*(\M^{\leq 3}_0) \rat \xrightarrow{\psi} \prod^{3}_{\delta=0} A^*(\M^{\delta}_0) \rat
$$
as the product of the maps $i^{3*}$, $i^{2*} j^{3*}$, $i^{1*} j^{2*}j^{3*}$ and $j^{1*}j^{2*}j^{3*}$.
Let $a$ be an element of $A^*(\M^{\leq 3}_0) \rat$ different from zero. If $\psi(a)$ is zero then it cannot be in the image of $i^3_*$ consequently $j^{3*}(a) \in A^*(\M^{\leq 2}_0)$ is different from zero. We can continue till we obtain that $j^{1*}j^{2*}j^{3*}$ is different from zero: absurd.
\end{dem}

\subsection{Restriction of classes to strata} \label{restrcla}
\medskip

Let us consider two trees $\Gamma$ and $\Gamma'$ with maximal multiplicity at most $3$ and number of edges respectively equal to $\delta$ and $\delta'$.

\begin{defi}
Given two graphs as above we call an ordered deformation of $\Gamma$ into $\Gamma'$ any surjective map of vertices $d: \Gamma' \to \Gamma$ such that
\begin{enumerate}
\item for each $P,Q \in \Gamma'$ we have $d(P) = d(Q)=A \in \Gamma'$ only if for each $R$ in the connected path from $P$ to $Q$ we have $d(R)=A$;
\item for each edge $(P,Q) \in \Gamma'$ such that $d(P) \neq d(Q)$ there must be an edge in $\Gamma$ between $d(P)$ and $d(Q)$.
\end{enumerate}
We denote by $\text{def}_o(\Gamma, \Gamma')$ the set of deformations.
\end{defi}

\begin{exm} \label{exmdef}
Let $\Gamma$ and $\Gamma'$ be the following graphs
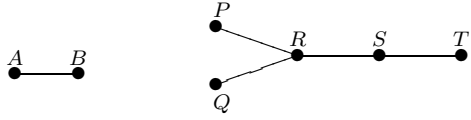
\begin{equation*}
\vcenter{\xymatrix { *=0{\bullet} \ar@{-}[r]^<A & *=0{\bullet} \ar@{}[r]^< B& *=0{}  }}
\qquad
\vcenter{\xymatrix @R=4pt {*=0{\bullet} \ar@{-}[dr]^<P & & & &\\
& *=0{\bullet} \ar@{-}[r]^<R &*=0{\bullet} \ar@{-}[r]^<S &*=0{\bullet} \ar@{}[r]^< T & *=0{}\\
*=0{\bullet} \ar@{-}[ur]_<Q & & & &}}
\end{equation*}

we have the following 8 ordered deformations

\begin{tabular}{|ll|ll|}
\hline
$1) (P,R,S,T) \mapsto A$ & $Q \mapsto B$ & $5) Q  \mapsto A$ & $(P,R,S,T)  \mapsto B$\\
\hline
$2) (Q,R,S,T) \mapsto A$ & $P \mapsto B$ & $6) P  \mapsto A$ & $(Q, R, S, T)  \mapsto B$\\
\hline
$3) (P,Q,R) \mapsto A$ & $(S,T) \mapsto B$ & $7) (S,T)  \mapsto A$ & $(P, Q, R)  \mapsto B$\\
\hline
$4) (P,Q,R,S) \mapsto A$ & $T \mapsto B$ & $8) T  \mapsto A$ & $(P, Q, R, S) \mapsto B$\\
\hline
\end{tabular}

\end{exm}

There exist two different equivalence relations in $\text{def}_o(\Gamma, \Gamma')$.
We say that two elements $d_1,d_2$ are in $\sim_\Gamma$ if there exists a $\gamma \in \aut(\Gamma)$ such that $d_2= \gamma d_1$. Similarly we say that two elements $d_1,d_2$ are in $\sim_{\Gamma'}$ if there exists a $\gamma' \in \aut(\Gamma')$ such that $d_2= d_1\gamma'$.

\begin{defi}
We call $\Gamma-$deformations (or simply deformations) from $\Gamma$ to $\Gamma'$ the set
$$
\text{def}_{\Gamma}(\Gamma, \Gamma') := \text{def}_o(\Gamma, \Gamma')/ \sim_{\Gamma}.
$$
We call $\Gamma'-$deformations from $\Gamma$ to $\Gamma'$ the set
$$
\text{def}_{\Gamma'}(\Gamma, \Gamma') := \text{def}_o(\Gamma, \Gamma')/ \sim_{\Gamma'}.
$$

\end{defi}

In the above example we can take as representatives of $\Gamma-$deformations the first 4 ordered deformations. On the other hand we have $1 \sim_{\Gamma'} 2$ and $5 \sim_{\Gamma'} 6$.  

Topological arguments let us state the following
\begin{propos}\label{defclaim}
Let $C \xrightarrow{\pi} T$ be a family of rational nodal curves over an irreducible scheme $T$. Suppose further that the generic fiber has  topological type $\Gamma$, then there exists a fiber of topological type $\widehat{\Gamma}$ only if there exists an ordered deformation of $\Gamma$ into $\widehat{\Gamma}$.
\end{propos}

Now let us consider the \'etale map
$$
\widetilde{\M}^{\Gamma}_0 \xrightarrow{\phi} \M^{\Gamma}_0.
$$
We have given above a description of $A^*{\M_0}$ as the subring of polynomials of $A^*(\widetilde{\M}^{\Gamma}_0)$ in the classes (corresponding to sections of $\widetilde{\M}_0^{\Gamma}$) $t_1, \dots, t_{\delta_1}, r_1, \dots, r_{\delta_2}$ invariant for the action of $\aut(\Gamma)$.

We call $\Ms_{0,i}$ the stack of rational nodal curves with $i$ sections.
Let $\left( \widetilde{\M}^{\Gamma}_0 \right)^{\leq \delta' -\delta}$ be the substack of
$$
\left( \Ms_{0,1} \right)^{\Delta_1} \times \left( \Ms_{0,2} \right)^{\Delta_2} \times \left( \Ms_{0,3}\right)^{\Delta_3}
$$
whose fibers have at most $\delta' -\delta$ nodes (the sum of nodes is taken over all the connected components).
Polynomials in $\Q[t_1, \dots, t_{\delta_1}, r_1, \dots, r_{\delta_2}]$ has a natural extension to $\left( \widetilde{\M}^{\Gamma}_0 \right)^{\leq \delta' -\delta}$. Let us fix one of such polynomials $a$ which are invariants for the action of $\aut{\Gamma}$.

The \'etale covering
$$
\widetilde{\M}^{\Gamma}_0 \xrightarrow{\phi} \M^{\Gamma}_0
$$
is obtained by gluing sections in a way which depends on $\Gamma$.

By gluing sections in the same way we obtain a functor
$$
\left( \widetilde{\M}^{\Gamma}_0 \right)^{\leq \delta' -\delta} \xrightarrow{\Pi} \M^{\delta'}_0.
$$

\begin{corol}\label{cordefclaim}
With the above notation we have that the closure $\overline \M^{\Gamma}_{0}$ of $\M^{\Gamma}_0$ in $\M^{\leq \delta'}_0$ is
$$
\Pi \left( \left( \widetilde{\M}^{\Gamma}_0 \right)^{\leq \delta' -\delta} \right).
$$
\end{corol}
\begin{dem}
From Proposition \ref{defclaim} we have
$$
\overline \M^{\Gamma}_0 \subseteq \Pi \left( \left( \widetilde{\M}^{\Gamma}_0 \right)^{\leq \delta' -\delta} \right).
$$
On the other hand let $C_{\Omega} \xrightarrow{\pi} \spk$ be the image in $\M^{\delta'}_0$ of a geometric point of $\left( \widetilde{\M}^{\Gamma}_0 \right)^{\leq \delta' -\delta}$. The dual graph $\Gamma_{\Omega}$ of $C_{\Omega}$ is a deformation of $\Gamma$. In order to show that $C_{\Omega} \xrightarrow{\pi} \spk$ is a geometric point of $\overline{\M}_0^{\Gamma}$, we fix a deformation $d: \Gamma_{\Omega} \to \Gamma$. For each vertex $A$ of $\Gamma$, the set $d^{-1}(A)$ is a subtree of $\Gamma_{\Omega}$. We can give a deformation $C_A \xrightarrow{\pi} T$ of $C_{\Omega}$ such that the generic fiber is $\Pro^1$. Furthermore we can define on $C_A \xrightarrow{\pi} T$ a family of $E(A)$ that respect $d$. At last we glue all $C_A$ along sections and obtain a deformation $C \xrightarrow{\pi} T$ of  $C_{\Omega} \xrightarrow{\pi} \spk$ in $\overline{\M}_0^{\Gamma}(\spk)$.
\end{dem}

\begin{defi}
We call the image of $\Pi$:
$$
\M_0^{\text{def}(\Gamma, \delta')}:=  \Pi \left( \left( \widetilde{\M}^{\Gamma}_0 \right)^{\leq \delta' -\delta} \right) \subset \M_0^{\leq \delta'}.
$$
\end{defi}

\begin{propos}
The map
$$
\Pi: \left( \widetilde{\M}_0^{\Gamma}\right)^{\delta'- \delta} \to \M_0^{\leq \delta'}
$$
is finite.
\end{propos}
\begin{dem}
({\bf sketch}) We have to prove that $\Pi$ is representable, with finite fibers and proper. The not trivial property to verify is properness. We can prove it through the valutative criterion (see \cite{hrt} p.101).
\end{dem}
Therefore the map $\Pi$ is finite hence projective, so we have the push-forward \cite{kre}
$$
\Pi_*: A^* \left( \left( \widetilde{\M}_0^{\Gamma}\right)^{\delta'- \delta} \right) \to A^* \left( \M_0^{\text{def}(\Gamma, \delta')}\right)
$$

\begin{defi}
Let $\Gamma$ be a tree with maximal multiplicity $\leq 3$ and $\Gamma'$ a deformation of $\Gamma$ (with maximal multiplicity $\leq 3$) with $\delta'$ edges.
Let $a$ be a class in $A^*(\M^{\Gamma}_0) \rat$ and $\tilde a$ its $\delta'$-lifting.
With reference to the cartesian diagram
\begin{equation*}
\xymatrix@=3pc{
\M^{\Gamma'}_0 \times_{\M^{\delta'}_0} \left ( \widetilde{\M}^{\Gamma}_0 \right )^ {\delta' -\delta} \cart \ar[r]^-{pr_2} \ar[d]_{pr_1} &\left ( \widetilde{\M}_0^{\Gamma} \right ) ^{\leq \delta' - \delta} \ar[d]^{\Pi}\\
\M_0^{\Gamma'} \ar[r]^{in} &\M_0^{\leq \delta'}
}
\end{equation*}
we define
$$
\Psi(\Gamma, \Gamma'):= \M^{\Gamma'}_0 \times_{\M^{\delta'}_0} \left ( \widetilde{\M}^{\Gamma}_0 \right )^ {\delta' -\delta}.
$$
\end{defi}

Topological arguments show the following
\begin{propos}\label{diagstrat}
$\Psi(\Gamma, \Gamma')$ is a disjoint union of components that we write as
$$
\Psi(\Gamma, \Gamma')= :\coprod_{\xi \in \text{def}_{\Gamma'}(\Gamma, \Gamma')} \Psi(\Gamma, \Gamma')_\xi
$$
where the union is taken over the set of ordered deformations up to $\sim_{\Gamma'}$.
\end{propos}

\begin{exm} $\;$

Let $\Gamma$ and $\Gamma'$ be the graphs of the Example \ref{exmdef}. We have $\delta'-\delta=3$ and
$$
\widetilde{\M}^{\Gamma}_0= \Ms^0_{0,1} \times \Ms^0_{0,1}
$$
with a double covering
$$
\widetilde{\M}^{\Gamma}_0 \xrightarrow{\phi} \M^{\Gamma}_0.
$$
Clearly $\Psi(\Gamma, \Gamma')$ is an inclusion of components
in 
\begin{equation}
\coprod _{i,j: i+j=3} \left( \Ms^i_{0,1} \times \Ms^j_{0,1} \right)
\end{equation}

We have that $\Psi(\Gamma, \Gamma')$ has 6 components which corresponds to deformations (enumerated in Example \ref{exmdef}) up to $\sim_{\Gamma'}$, with the following inclusions:
\begin{itemize}
\item deformation 1 (which is $\Gamma'$-equivalent to 2) and 4 correspond to two connected components of $\Ms^{3}_{0,1} \times \Ms^{0}_{0,1}$
\item deformation 3 corresponds to a connected component of $\Ms^{2}_{0,1} \times \Ms^{1}_{0,1}$
\item deformation 7 corresponds to a connected components of $\Ms^{1}_{0,1} \times \Ms^{2}_{0,1}$
\item deformation 5 (which is $\Gamma'$-equivalent to 6) and 8 corresponds to two connected components of $\Ms^{0}_{0,1} \times \Ms^{3}_{0,1}$
\end{itemize}
\end{exm}

\begin{defi}
Let us consider a class $a$ in $A^*(\M^{\Gamma}_0) \rat$.
With reference to the \'etale covering
$$
\widetilde{\M}^{\Gamma}_0 \xrightarrow{\phi} \M^{\Gamma}_0
$$
we define the {\it lifting} of $a$
$$
\widetilde{a}:= \frac{\phi^*a}{\sigma}\in A^*(\widetilde{\M}^{\Gamma}_0) \rat,
$$
this means that $\phi_*(\widetilde{a})=a$.
Since $\widetilde a$ can be written as a polynomial in the classes $\psi$ defined in Proposition \ref{ringstratum}
that depends only on $\Gamma$, it has a natural extension to $A^*\left( \left( \widetilde{\M}_0^{\Gamma}\right)^{\delta'- \delta} \right) \rat$ that we call $\delta'$-lifting of $a$ and we still write $\tilde a$.
\end{defi}

For every $\xi \in \text{def}_{\Gamma'}(\Gamma, \Gamma')$ we have the related commutative diagram
\begin{equation*}
\xymatrix@=3pc{
\Psi(\Gamma, \Gamma')_\xi \ar[r]^{pr_2^\xi} \ar[d]_{pr_1^\xi} & \left ( \widetilde{\M}_0^{\Gamma} \right ) ^{\leq \delta' - \delta} \ar[d]^{\Pi}\\
\M_0^{\Gamma'} \ar[r]^{in} &\M_0^{\leq \delta'}
}
\end{equation*}
The map $pr_2^\xi$ is a closed immersion of codimension $\delta' - \delta$. Let us still call $\widetilde a$ the pullback of the polynomial $\widetilde a$ through $pr_2^{\xi}$.
By the excess intersection formula (Section 6.3 \cite[ful], similar arguments show it for algebraic stacks)
we have
$$
in^!(\widetilde a)= \sum_{\xi \in \text{def}_{\Gamma'}(\Gamma, \Gamma')} (\widetilde a \cdot c_{top}[\Nf^{\xi}]).
$$
where $\Nf^{\xi}:=(pr_1^{\xi*} \Nf_{in})/ \Nf_{pr_2^{\xi}}$.

For each $\xi \in \text{def}_{\Gamma'}(\Gamma, \Gamma')$ we have an \'etale covering
$$
f_{\xi}: \widetilde \M^{\Gamma'} \to \Psi(\Gamma,\Gamma')_{\xi}
$$
that glue along sections.

For every $\xi \in \text{def}_{\Gamma'}(\Gamma, \Gamma')$ let us call $\sigma'_{\xi}$ the degree of $pr^{\xi}_1$.

We have the following commutative diagram
\begin{equation*}
\xymatrix{
\widetilde{\M}_0^{\Gamma'} \ar[rr]^{f_{\xi}} \ar[dr]_{\phi'} && \Psi(\Gamma, \Gamma')_\xi \ar[dl]^{pr_1^\xi}\\
& \M_0^{\Gamma'} &
}
\end{equation*}
from which we have
$$
\text{ord} \phi' = \text{ord} f_{\xi} \cdot \text{ord} \sigma'_{\xi}.
$$
\begin{oss} \label{restrpol}
The order of $f_{\xi}$ is the number of $g \in \aut(\Gamma')$ such that for each deformation $d$ associated to $\xi$ the deformation $d \circ g$ is $\Gamma'$-equivalent to $d$.
\end{oss}
We want to explicit ${\phi'}^*pr_{1*}in^!(\widetilde a)$ in $A^* ( \widetilde{\M}_0^{\Gamma'} ) \rat$: being invariant for the action of $\aut(\Gamma')$ we can see it as a class in $\M^{\Gamma'}_0$ which is the restriction of the extension of the class $a \in A^* \left( \M^{\Gamma}_0 \right) \rat$.

\begin{propos}\label{restriction}
We have
\begin{equation*}
{\phi'}^*pr_{1*}in^!(\widetilde a)=\sum_{\xi \in \text{def}_{\Gamma'}(\Gamma, \Gamma')} \sigma'_{\xi} f^*_{{\xi}}(\widetilde a) \cdot c_{top}(\widetilde{\Nf}^{\xi})
\end{equation*}
where
$$
\widetilde{\Nf}^{\xi}:= f^*_{{\xi}} {\Nf}^{\xi}.
$$
\end{propos}

\begin{dem}
Putting together all the above remarks and definitions we have
\begin{eqnarray*}
{\phi'}^*pr_{1*}in^!(\widetilde a)&=& {\phi'}^*pr_{1*} \sum_{\xi \in \text{def}_{\Gamma'}(\Gamma, \Gamma')} (\widetilde a \cdot c_{top}[\Nf^{\xi}])\\
&=&  {\phi'}^* \sum_{\xi \in \text{def}_{\Gamma'}(\Gamma, \Gamma')} pr^{\xi}_{1*}(\widetilde a \cdot c_{top}[\Nf^{\xi}])\\
&=& \sum_{\xi \in \text{def}_{\Gamma'}(\Gamma, \Gamma')} f^{*}_{\xi} pr^{\xi*}_{1} pr^{\xi}_{1*}(\widetilde a \cdot c_{top}[\Nf^{\xi}])\\
&=& \sum_{\xi \in \text{def}_{\Gamma'}(\Gamma, \Gamma')} \sigma'_{\xi} f^*_{{\xi}}(\widetilde a) \cdot c_{top}(\widetilde{\Nf}^{\xi})
\end{eqnarray*}
\end{dem}

We explicit now the computation of classes corresponding to strata.

\begin{propos}\label{refrestr}
Given a tree $\Gamma$ with maximal multiplicity $\leq 3$, let $\gamma_{\Gamma}$ be the class of $\overline \M^{\Gamma}_0$ in $\M_0$ and let $\Gamma'$ be another tree with maximal multiplicity $\leq 3$.

If $\Gamma'$ is a deformation of $\Gamma$, then the restriction of $\gamma_{\Gamma}$ to $A^* \left( \M^{\Gamma'} \right) \rat$ is (following the above notation)
\begin{equation*}
\sum_{\xi \in \text{def}_{\Gamma'}(\Gamma, \Gamma')} \frac{\sigma'_{\xi}}{\sigma}  c_{top}(\widetilde{\Nf}^{\xi}).
\end{equation*}

If $\Gamma'$ is not a deformation of $\Gamma$ then the restriction is 0.
\end{propos}

\begin{dem}
In case $\Gamma'$ is a deformation of $\Gamma$, the polynomial $a$ of Proposition (\ref{restriction}) is 1, consequently $\widetilde a$ is $1/\sigma$. With reference to the \'etale covering $\widetilde{\M}^{\Gamma'}_0 \xrightarrow{\phi} \M^{\Gamma'}_0$ and the inclusion $\M^{\Gamma'}_0 \xrightarrow{in} \M^{\delta'}_0$, from Proposition \ref{restriction} we obtain
\begin{equation*}
\sum_{\xi \in \text{def}_{\Gamma'}(\Gamma, \Gamma')} \frac{\sigma'_{\xi}}{\sigma}  c_{top}(\widetilde{\Nf}^{\xi}).
\end{equation*} 
If $\Gamma'$ is not a deformation of $\Gamma$ then using Corollary \ref{cordefclaim} and basic topological arguments, we have that $\M^{\Gamma'}_0$ does not intersect the closure of $\M^{\Gamma}_0$ in $\M_0$, so the restriction of the class must be 0.
\end{dem}

\begin{oss}
For each $\xi \in \text{def}_{\Gamma'}(\Gamma, \Gamma')$, we have an exact sequence of sheaves
$$
0 \to f^*_{\xi} \Nf_{pr^{\xi}_2} \to  {\phi'}^* \Nf_{in} \to \widetilde \Nf^{\xi} \to 0.
$$
We have seen in Proposition \ref{classstratum} how to compute $c_{top} ({\phi'}^* \Nf_{in})$, similarly we can compute $c_{top} (f^*_{\xi} \Nf_{pr^{\xi}_2})$ and finally we consider the following relation (that follows from the exact sequence)
$$
c_{top} ({\phi'}^* \Nf_{in})= c_{top} (f^*_{\xi} \Nf_{pr^{\xi}_2}) \cdot c_{top} (\widetilde \Nf^{\xi}).
$$
\end{oss}
We carry on the calculation for trees with at most three nodes in the last Section.

\subsection{Mumford classes}\label{mumcla}

Given a tree $\Gamma$ with at most four vertices, we have the restriction of the dualizing sheaf $\omega_0 := \omega_{\cu/ \M_0}$ on the universal curve $\cu^{\Gamma}$ of $\M^{\Gamma}_0$ (see Section \ref{normal}). We  call $\omega^{\Gamma}_0:=\omega_{\cu^{\Gamma} / \M^{\Gamma}_0}$ its restriction.

We consider two kinds of classes on $A^*_{\aut(C_0)}$ induced by the sheaf $\omega^{\Gamma}_0$:
\begin{enumerate}
\item the pushforward of polynomials of the Chern class $K:=c_1(\omega^{\Gamma}_0)$;
\item the Chern classes of the pushforward of $\omega_0$.
\end{enumerate}  

The first kind will give us the equivalent of Mumford classes, but on $\M_0^{\leq 3}$ we can only define classes of the second type (see Section \ref{normal}).

The aim of this section is to compute classes of the first kind for $\Gamma$ with at most three nodes  and to describe them as polynomials in classes of the second kind. If this description is independent from the graph $\Gamma$ then we can define them as elements of $A^*(\M_0^{\leq 3})$.

First of all we define on $\cu^{\Gamma}$ and $\M^{\Gamma}_0$ the following classes
\begin{eqnarray*}
K &:=& c_1(\omega^{\Gamma}_0) \in A^1(\cu^{\Gamma}),\\
\km_i &:=& \Pi_*(K^{i+1}) \in A^{i}(\M^{\Gamma}_0).
\end{eqnarray*}

Such classes $\km_i$ (introduced in \cite{mumenum} for the moduli spaces of stable curves) are called {\it Mumford classes}.

{\bf Mumford classes on $\M_0^0$}

The stack $\M^0_0$ is $\B \Pro Gl_2$.

Let us consider the universal curve
$$
[\Pro ^1 _{\C} / \Pro Gl_2] \xrightarrow{\Pi} \B \Pro Gl_2.
$$
We call $\overline \Pi$ the induced map $\Pro ^1 _{\C} \to \spc$ and $\overline \omega^0_0$  a lifting of $\omega^0_0$ on $\Pro^1_{\C}$.

We notice that $\left ( \overline \omega^0_0 \right ) ^{\vee} = T _ {\overline \Pi}$ ($T_{\Pi}$ is the relative tangent bundle along $\Pi$), we set  $$ K:=c^{\Pro Gl_2}_1(T _{\overline \Pi})= -c^{\Pro Gl_2}_1(\overline \omega^0_0).$$

Furthermore  $\overline \Pi_*(T_{\overline \Pi})=H^0(\Pro^1_{\C}, T_{\overline \Pi})=\sla _2$ seen as adjoint representation of $\Pro Gl_2$.

By applying the equivariant Grothendieck-Riemann-Roch Theorem we obtain
\begin{eqnarray*}
ch(\Pi_*(\omega_0^\vee)) &=& \Pi_*(Td(T_{\Pi})ch(\omega_0^\vee))\\
ch^{\Pro Gl_2}(\sla_2)&=&\Pi_*(Td^{\Pro Gl_2}(T_{\overline \Pi})ch^{\Pro Gl_2}(T_{\overline \Pi}))\\
3 -c^{\Pro Gl_2}_2(\sla _2)&=& \Pi_* \left[ \left( e^{- (K)}\right) \left( \frac{-K}{1- e^{K}}\right) \right]
\end{eqnarray*}
By applying GRR to the trivial linear bundle we obtain:
$$
1= \Pi_* \left[ \frac{-K}{1- e^{K}} \right]
$$
If we subtract the second equation from the first, we obtain
$$
2 -c^{\Pro Gl_2}_2(\sla _2)  = \Pi_* \left[ \left( \frac{1-e^{K}}{e^{K}}\right) \left( \frac{-K}{1- e^{K}}\right) \right] = \Pi_* \left[ -K e^{-K} \right],
$$
from which we get the following
\begin{propos}\label{mumzero}
On $\M^0_0$ we have
\begin{eqnarray*}
\km_0 = -2, \quad \km_2 = 2c^{\Pro Gl_2}_2(\sla _2), \quad \km_1 = \km_3 = 0.
\end{eqnarray*}
\end{propos}

{\bf Mumford classes on strata of singular curves}

Now let us consider the following cartesian diagram (see Section \ref{stratclass})
\begin{equation}\label{diagtildef}
\xymatrix@=3pc{
[\stackrel{\delta + 1}{\coprod} \Pro ^1 _{\C} /  H] \cart \ar@/_3pc/[dd]_{\widetilde F} \ar [r] ^\xi \ar[d]_{\widetilde N} & \widehat{\cu} ^{\Gamma} \ar@/^2pc/[dd]^{F} \ar[d]^{N} \\
[C_0 / H]  \cart \ar [d]_{\widetilde \Pi} \ar[r] &\cu ^{\Gamma} \ar[d]^\Pi\\
\B H \ar[r] _{\phi} & \B  \aut(C_0)
}
\end{equation}
where $N: \widehat{\cu} ^{\Gamma} \to \cu^{\Gamma}$ is the normalization of $\cu^{\Gamma}$ described in Section \ref{normal}.

In the following we call $\widetilde{\omega}^{\Gamma}_0$ the sheaf
$$
\xi^* N^*(\omega^{\Gamma}_0).
$$ 
and $\widetilde{K}:=c_1(\omega_0)$.

\begin{propos} \label{pfkm}
Using the above notation, the Mumford classes $\km_i \in A^i(\B F \ltimes H) \rat$ are described by the following relation
$$
\phi_* \widetilde{F}_*(\widetilde{K}^{i+1}) = \sigma \km_i,
$$
\end{propos}

\begin{dem}
Let us fix an index $i \in \N$. Since the map $N$ is finite and generically of degree 1 we have $N_* N^* = \id$
therefore $\km_i := \Pi_* K^{i+1} = \Pi_*(N_* N^*)K^{i+1} = F_*(N^* K^{i+1})$,
since $F$ is projective and $\phi$ is \'etale we can apply the projection formula and obtain $$\phi^* F_*(N^* K^{i+1})= \widetilde{F}_* \xi ^* (N^* K^{i+1}).$$
Furthermore the map $N \circ \xi$ is finite and so $$\xi^* N^* c_1(\omega^{\Gamma}_0)^{i+1} = c_1(\xi^* N^* (\omega^{\Gamma}_0))^{i+1} = \widetilde K ^{i+1}.$$
We conclude by noting that $\phi_* \phi^*$ is multiplication by $\sigma$.
\end{dem}

\begin{propos}\label{defmumstrat}
Following notation of Proposition \ref{classstratum} and order elements of $\Delta_1$ from $1$ to $\delta_1$, we have
\begin{eqnarray} \label{kmum}
\km_m=-\phi_*\frac{t_1^m+ \dots + t_{\delta_1}^m}{\sigma}.
\end{eqnarray}
\end{propos}

\begin{dem}
From Proposition \ref{pfkm}, we have reduced the problem to computing $\widetilde K$ and then writing pushforward along $\widetilde{F}$. With respect to each component of the stack
$$
\widehat{\cu} ^{\Gamma} = \stackrel{\delta + 1}{\coprod} [\Pro ^1 _{\C} / \aut(C_0)]
$$
we fix coordinates on $\Pro ^1 _{\C}$ which are compatible with coordinates chosen on $C$. By forgetting the action of $\aut(\Gamma)$ we keep the same system of coordinates on $\stackrel{\delta + 1}{\coprod} \Pro ^1 _{\C}$ when we consider $\stackrel{\delta + 1}{\coprod} [\Pro ^1 _{\C} / H]$.

Giving a bundle on a quotient stack $[X/G]$ is equivalent to giving a bundle  $U \to X$ equivariant for the action of $G$. 
By abuse of notation we still call $\widetilde{\omega}^{\Gamma}_0$ any lifting of $\widetilde{\omega}^{\Gamma}_0$ on $\stackrel{\delta + 1}{\coprod} \Pro ^1 _{\C}$.

In order to make computations we need to render explicit the action of $H:=E^{\Delta_1} \times \Gm^{\Delta_2}$ on
$$
\widetilde{F}_*(\widetilde{\omega}_0)= H^0 \left( \stackrel{\delta + 1}{\coprod} \Pro ^1 _{\C}, \widetilde{\omega}^{\Gamma}_0 \right ).
$$
 
Set $\Delta:=\Delta_1 \cup \Delta_2$, from the inclusion $\Gm \to E$ we have a cartesian diagram of stacks
\begin{equation*}
\xymatrix{
\B (\Gm)^{\Delta} \cart \ar[r]^-{\phi} \ar[d]^{\Psi} &\B (\aut(\Gamma) \ltimes (\Gm)^{\Delta}) \ar[d]^{\Psi}\\
\B H \ar[r]^{\phi} &\B \aut (C_0)
}
\end{equation*}
Roughly speaking we can say that we obtain the stacks in the top row by fixing the point $0$ on components with a node. The functor $\Psi$ forgets these points.
We have the following ring isomorphisms 
\begin{eqnarray*}
A^*_{H} &\xrightarrow{\Psi^*}& A^*_{(\Gm)^{\Delta}} \\
 A^*_{\aut(C_0)} &\xrightarrow{\Psi^*}& A^*_{\aut(\Gamma) \ltimes (\Gm)^{\Delta}}.
\end{eqnarray*}
We have defined the classes $t_{\Delta_1}, r_{\Delta_2}$ in $A^*_H$. By using the same notation of Section \ref{stratclass}, the map $\Psi^*$ is the identity on $r_{\Delta_2}$. For each vertex $P$ in $\Gamma$ such that $e(P)=1$ the map $\Psi^*$ sends $t_P$ to $c^{\Gm}_1(H^0(\Pro^1_{\C},\Of(z_\infty)))$ of the same component, which we still call $t_P$.
Consequently, with reference to the following diagram
\begin{equation*}
\xymatrix@=3pc{
[\stackrel{\Delta}{\coprod} \Pro ^1 _{\C} /  (\Gm)^{\Delta}] \ar[d]_{\widetilde{F}} \\
\B (\Gm)^{\Delta} \ar[r] _-{\phi} & \B  (\aut(\Gamma) \ltimes (\Gm)^{\Delta})
}
\end{equation*}
we have reduced the problem to consider the action of $(\Gm)^{\Delta}$ on
$$
\widetilde{F}_*(\widetilde{\omega}^{\Gamma}_0)= H^0 \left( \stackrel{\Delta}{\coprod} \Pro ^1 _{\C}, \widetilde{\omega}^{\Gamma}_0 \right )
$$
where, again with abuse of notation, we call $\widetilde{\omega}^{\Gamma}_0$ the sheaf $\Psi^*(\widetilde{\omega}^{\Gamma}_0)$.

The map $\widetilde{F}$ is the union of maps	
$$
\widetilde{F}_P:= [\Pro^1 _{\C}/ (\Gm)^\Delta] \to \B (\Gm)^\Delta
$$
where only the component of $(\Gm)^{\Delta}$ corresponding to $P \in \Delta$ does not acts trivially on $\Pro^1_{\C}$ and the action is
\begin{eqnarray*}
\az: \Gm \times \Pro^1_{\C} &\to& \Pro^1_{\C}\\
(\lambda, [X_0,X_1]) &\mapsto& [X_0, \lambda X_1].
\end{eqnarray*}
Set
$$
P_0:=[0,1], \quad P_1:=[1,1], \quad P_{\infty}:= [1,0].
$$
We have that $\mathcal F := (\widetilde{\omega}_0)^\vee$ is the sheaf on  $\stackrel{\Delta}{\coprod}\Pro ^1 _{\C}$ such that restricted:
\begin{itemize}
\item to the components with one node it is $\mathcal F_P:=(\omega \otimes \Of (z_{\infty}))^\vee$\\
\item to the components with two nodes it is $\mathcal F_P:=(\omega \otimes \Of (z_{\infty} + z_0))^\vee$\\
\item to the components with three nodes it is $\mathcal F_P:=(\omega \otimes \Of (z_{\infty} + z_0 + z_1))^\vee$
\end{itemize}
where $\omega$ is the canonical bundle on $\Pro^1_{\C}$.

In the following, Chern classes will be equivariant for the action of $\Gm$. For each $P \in \Delta$ let us indicate with $K_P$ the class $c^{\Gm}_1(\omega) \in A^1_{\Gm}(\Pro^1_{\C})$ on each component, with $R_P$ the class $c^{\Gm}_1(\Of(z_{\infty})) \in A^1_{\Gm}(\Pro^1_{\C})$, and with $Q_P$ the class $c^{\Gm}_1(\Of(P_0)) \in A^1_{\Gm}(\Pro^1_{\C})$.

We have
\begin{eqnarray*}
c^{\Gm}_1(\F_{\Delta_1})&=& - {K_{\Delta_1}} - R_{\Delta_1}\\
c^{\Gm}_1(\F_{\Delta_2})&=& - {K_{\Delta_2}} - R_{\Delta_2} - Q_{\Delta_2}.
\end{eqnarray*}
In order to determine Mumford classes it is then necessary to compute the pushforward classes
\begin{eqnarray*}
&& \widetilde{F}_{\Delta_1*} (- {K_{\Delta_1}} - R_{\Delta_1})^h \text{ and}\\
&& \widetilde{F}_{\Delta_2*} (- {K_{\Delta_2}} - R_{\Delta_2} - Q_{\Delta_2})^h
\end{eqnarray*}
for every natural $h$.

Let us start with computing the push-forward along $\widetilde F _{\Delta}$ of every power of $K_{\Delta}=c^{\Gm}_1(\omega )$.

We have
$$
\widetilde F_{\Delta_1*}(\omega^\vee)=H^0(\Pro ^1 _\C, \omega^\vee).
$$
Now it is necessary to determine the action of $(\Gm)$
on global sections of  $\omega ^\vee$. Let $z:= X_1 / X_0$ be the local coordinate around $z_{\infty}$, the global sections of $\omega^{\vee}$ are generated as a vectorial space by $\hol, z \hol, z^2 \hol$. Relatively to this basis, the action of $\Gm$ is given by
\begin{eqnarray*}
\az: \Gm \times H^0(\Pro^1_{\C},\omega^\vee) &\to& H^0(\Pro^1_{\C}, \omega^\vee)\\
(\lambda, a\hol + b z \hol + c z^2 \hol) &\mapsto& (a \lambda \hol + b z \hol + c \frac{1}{\lambda} z^2 \hol)
\end{eqnarray*}
the multiplicity of the action is therefore $(1,0,-1)$, so the Chern character of $\widetilde F _{\Delta*}(\omega^\vee)$ is $e^{t_{\Delta}} + 1 + e^{-t_{\Delta}}$.

Let us observe that $\omega^{\vee}$ is the tangent bundle relative to $\widetilde F _{\Delta}$, so by applying the GRR Theorem we get
$$
1 + e^{t_{\Delta}} + e^{-t_{\Delta}} = \widetilde F_{\Delta*} \left[ \left( e^{- (K_{\Delta})}\right) \left( \frac{-K_\Delta}{1- e^{K_\Delta}}\right) \right]
$$
By applying GRR to the trivial linear bundle we obtain:
$$
1= \widetilde F_{\Delta*} \left[ \frac{-K_\Delta}{1- e^{K_\Delta}} \right]
$$
If we subtract the second equation from the first, we obtain
$$
e^{t_\Delta} + e^{-t_\Delta}  = \widetilde F_{\Delta*} \left[ \left( \frac{1-e^{K_\Delta}}{e^{K_\Delta}}\right) \left( \frac{-K_\Delta}{1- e^{K_\Delta}}\right) \right] = \widetilde F_{\Delta*} \left[ -K_\Delta e^{-K_\Delta} \right],
$$
from this, by distinguishing between even and odd cases, it follows that

\begin{eqnarray*}
\widetilde F_{\Delta*}K_\Delta^{2h}&=&0\\
\widetilde F_{\Delta*}K_\Delta^{2h+1}&=& -2t_\Delta^{2h}
\end{eqnarray*}

Let us notice that there exist two equivariant sections of $\widetilde F_{\Delta}$ given by the fixed points $z_0$ and $z_\infty$, that we will call respectively $s_0$ and $s_\infty$
\begin{equation*}
\xymatrix{
\B (\Gm)^{\Delta} \ar@<0.5ex>[r]^-{s_0} \ar@<-0.5ex>[r]_-{s_{\infty}} & [\Pro^1_{\C}/ (\Gm)^{\Delta} ]
}
\end{equation*}
From the self intersection formula we have 
\begin{eqnarray*}
&&s_0^*(Q_\Delta)=-t_\Delta,  \; s_0^*(R_\Delta)=0, \; s_0^*(K_\Delta)=t_\Delta\\
&&s_\infty^*(Q_\Delta)=0, \; s_\infty^*(R_\Delta)=t_\Delta, \; s_\infty^*(K_\Delta)=-t_\Delta;
\end{eqnarray*}
then by applying the projection formula (see \cite{ful} p.34)
for every cycle $D \in A^*_{(\Gm)^\Delta}(\Pro^1)$ we have
\begin{eqnarray*}
&& Q_\Delta \cdot D = s_{0*}(1) \cdot D = s_{0*}(s_0^*D)\\
&& R_\Delta \cdot D = s_{\infty*}(1) \cdot D = s_{\infty*}(s_\infty^*D).
\end{eqnarray*} 
With reference to the previous relations, the Mumford classes are determined on every component (we separate between even and odd cases)
\begin{eqnarray*}
\widetilde F_{\Delta *}(-K_\Delta -R_\Delta)^{2h} &=& \widetilde F^\Delta_*(K_\Delta)^{2h} + 
					\sum_{a=1}^{2h} \binom{2h}{a}				
					\widetilde F_{\Delta*}(R_\Delta^{a}\cdot K_\Delta^{2h-a})\\
				&=& 	\sum_{a=1}^{2h} \binom{2h}{a}	
					(\widetilde F_{\Delta*}s_{\infty*})s_\infty^*(R_\Delta^{a-1}\cdot K_\Delta^{2h-a})\\
				&=& 		\sum_{a=1}^{2h} \binom{2h}{a}	
						(-1)^{a-1}t_\Delta^{2h-1} = -t_\Delta^{2h-1}
\end{eqnarray*}
and in a completely analogous way, we have the following relations
\begin{eqnarray*}
&&\widetilde F_{\Delta_1*}(-K_{\Delta_1} -R_{\Delta_1})^{2h + 1} = t_i^{2h}\\
&&\widetilde F_{\Delta_2*}(-K_{\Delta_2} -R_{\Delta_2} -Q_{\Delta_2})^{h}=0.
\end{eqnarray*}
The last relation allows us to ignore components with two nodes.
By following the notation of Proposition \ref{pfkm} we notice that
$$
\widetilde F _* \widetilde K ^ {m+1} = \sum_{P \in \Delta_1} \widetilde F_{P*}(-1)^{m+1}(-K_P -R_P)^{m+1} 
$$
and consequently, if we order elements of $\Delta_1$ from $1$ to $\delta_1$, we have

\begin{eqnarray} \label{kmum}
\km_m=-\phi_*\frac{t_1^m+ \dots + t_{\delta_1}^m}{\sigma}.
\end{eqnarray}

\end{dem}

{\bf Definition of Mumford classes on $\M_0^{\leq 3}$}
Now we consider the universal curve
$$
\cu^{\leq 3} \xrightarrow{\Pi} \M_0^{\leq 3}.
$$
We have seen in Section \ref{normal} that the pushforward $\Pi_* \left ( \omega^{\leq 3}_0 \right )^{\vee}$ is a well defined rank three vector bundle. Consequently from \cite{kre} Section 3.6 we have that $c_i(\Pi_* \left ( \omega^{\leq 3}_0 \right )^{\vee})=0$ for $i>3$.
We fix the following notation
\begin{eqnarray*}
 \cm_1:=  c_1(\Pi_*\left ( \omega^{\leq 3}_0 \right )^{\vee}), \quad \cm_2:= c_2(\Pi_* \left ( \omega^{\leq 3}_0 \right )^{\vee}), \quad \cm_3:= c_3(\Pi_* \left ( \omega^{\leq 3}_0 \right )^{\vee}).
\end{eqnarray*}
We still call $\cm_1, \cm_2, \cm_3$ their restriction to each stratum of $\M^{\leq 3}_0$.

\begin{defi}\label{defmumpol}
We define in $A^*(\M^{\leq 3}_0)$ Mumford classes $\km_1, \km_2, \km_3$ as follows
\begin{eqnarray*}
\km_1 :=  -\cm_1, \quad \km_2  :=  2 \cm_2 - \cm^2_1, \quad \km_3 := -\cm^3_1 + 3 \cm_1 \cm_2 -3\cm_3.
\end{eqnarray*}
\end{defi}

\begin{oss}
From Proposition \ref{injleqtre}, the Chow ring $A^*(\M^{\leq 3}_0) \rat$ injects into the product of $A^*(\M^{\Gamma}_0)\rat$ over trees with at most three edges. Consequently in order to verify that the above is a good definition we only need to prove that the restrictions of Mumford classes to each stratum are the given polynomials.
\end{oss}

\begin{propos}\label{finmum}
Let $\Gamma$ be a tree with at most three edges. Set 
\begin{eqnarray*}
 \cm_1(\Gamma):=  c_1(\Pi_*\left ( \omega^{\Gamma}_0 \right )^{\vee}), \; \cm_2(\Gamma):= c_2(\Pi_* \left ( \omega^{\Gamma}_0 \right )^{\vee}), \; \cm_3(\Gamma):= c_3(\Pi_* \left ( \omega^{\Gamma}_0 \right )^{\vee}).
\end{eqnarray*}
and
\begin{eqnarray*}
&& \nm_1(\Gamma) :=  \cm_1(\Gamma), \; \nm_2(\Gamma)  :=  \cm^2_1(\Gamma) -2 \cm_2(\Gamma), \\
&& \nm_3(\Gamma) := \cm^3_1(\Gamma) - 3 \cm_1(\Gamma) \cm_2(\Gamma) +3\cm_3(\Gamma).
\end{eqnarray*}

Then we have
\begin{eqnarray*}
\km_1 =  -\nm_1(\Gamma) \quad \km_2  =  -\nm_2(\Gamma) \quad \km_3 =  -\nm_3(\Gamma)
\end{eqnarray*}
\end{propos}

\begin{dem}
We consider first of all the case $\M^{\Gamma}_0= \M^0_0$. Here we have
\begin{eqnarray*}
&&\nm_1(\Gamma) = c^{\Pro Gl_2}_1(\sla_2)=0\\
&&\nm_2(\Gamma) =-2c^{\Pro Gl_2}_2(\sla_2)\\
&&\nm_3(\Gamma)=0
\end{eqnarray*}
and we can conclude because on $\M^0_0$ we have that $\km_1=\km_3=0$ and $\km_2=2c^{\Pro Gl_2}_2(\sla_2)$ (Proposition \ref{mumzero}).

Now let us consider any other tree $\Gamma$ with at most three edges.  With reference to diagram (\ref{diagtildef}),
as Chern classes commutes with base changing, we have
$$
\phi^*c_i(\Pi_*\left( \omega^{\Gamma}_0 \right)^{\vee}) = c_i(\widetilde \Pi_* \left( \widetilde \omega^{\Gamma}_0 \right) ^{\vee}).
$$
Given $\Gamma$ we order elements of $\Delta_1$ from $1$ to $\delta_1$. By putting together the above relation and the equation (\ref{kmum}), we reduce to show
$$
ch(\widetilde \Pi_* \left( \widetilde \omega^{\Gamma}_0 \right)^{\vee})= 3 + \sum^{\infty}_{m=1} \frac{t^m_1 + \dots + t^m_{\delta_1}}{m!}.
$$
We recall that the universal curve $\widetilde \cu^{\Gamma}$ on $\widetilde \M_0^{\Gamma}= \B (E^{\Delta_1} \times \Gm^{\Delta_2})$ is the quotient stack $[C_0/ E^{\Delta_1} \times \Gm^{\Delta_2}]$.

Since we have a morphism $\Psi: \B (\Gm)^{\Delta} \to \widetilde \M_0^{\Gamma}$
such that $\Psi^*: A^*(\widetilde \M_0^{\Gamma}) \to A^*_{(\Gm)^{\Delta}}$ is an isomorphism, with reference to the cartesian diagram
\begin{equation*}
\xymatrix{
[C_0 / \Gm^{\Delta}] \cart \ar[r]^-\Psi \ar[d]_{\widetilde \Pi} & [C_0/E^{\Delta_1} \times \Gm^{\Delta_2}] \ar[d]^{\widehat \Pi}\\
\B \Gm^{\Delta} \ar[r]^-\Psi & \B E^{\Delta_1} \times \Gm^{\Delta_2}
}
\end{equation*}
we reduce to consider the pullback sheaf $\Psi^* \left( \widetilde \omega^{\Gamma}_0 \right)^\vee$ which we still call  $\left( \widetilde \omega^{\Gamma}_0 \right)^\vee$ as its lifting to $C_0$.
As $H^1(C_0, \left( \widetilde \omega^{\Gamma}_0 \right)^\vee)=0$ (see proof of Preposition \ref{dualsh}), we have
\begin{equation*}
ch(\widetilde \Pi_* \left( \widetilde \omega^{\Gamma}_0 \right)^{\vee})= ch^{\Gm^\Delta}(H^0(C_0, \left( \widetilde \omega^{\Gamma}_0 \right)^\vee))
\end{equation*}
On curves of topological type
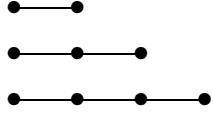
\begin{eqnarray*}
&& \vcenter{\unoa}\\
&& \vcenter{\duea}\\
&& \vcenter{\trea}
\end{eqnarray*}
we have that global sections of $ \left( \widetilde \omega^{\Gamma}_0 \right)^{\vee}$ are sections of $\Of_{\Pro^1_{\C}}(1)$ on extremal components and $\Of_{\Pro^1_{\C}}$ on the other components, which agree on nodes. Since $H^0(\Pro^1_\C,\Of_{\Pro^1_{\C}}(0))= \C$,  global sections of $ \left( \widetilde \omega^{\Gamma}_0 \right)^{\vee}$ are sections of $\Of_{\Pro^1_{\C}}(1)$ on the two extremal components which are equal on nodes. On each extremal component we fix coordinates $[X'_0,X'_1]$ and $[X''_0, X''_1]$. Sections on $\Of_{\Pro^1_{\C}}(1)$ are linear forms
\begin{eqnarray*}
&& a_1 X'_0 + b_1 X'_1\\
&& a_2 X''_0 + b_2 X''_1
\end{eqnarray*}
which agree at $z_{\infty}=[1,0]$. This happens if and only if  $a_1=a_2$. We have only $(\Gm)^{\Delta_1}$ which does not acts trivially on $H^0(C_0, \left( \widetilde \omega^{\Gamma}_0 \right)^{\vee})$ and the action is
\begin{eqnarray*}
\az: (\Gm \times \Gm) \times H^0(C_0, \left( \widetilde \omega^{\Gamma}_0 \right)^{\vee}) &\to& H^0(C_0, \left( \widetilde \omega^{\Gamma}_0 \right)^{\vee})\\
(\lambda_1,\lambda_2), (a, b_1, b_2 ) &\mapsto& ( a, \lambda_1 b_1, \lambda_2 b_2)
\end{eqnarray*}
so the Chern character of $\widetilde \Pi_* \omega^{\vee}_0$ is
$$
1 + e^{t_1} + e^{t_2} = 3 + \sum^{\infty}_{m=1} \frac{t^m_1+ t^m_2}{m!}
$$
as we claimed.

The last case to consider is when $\Gamma$ equals to
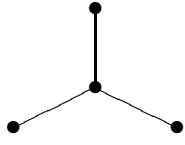
\begin{eqnarray*}
&& \vcenter{\treb}\\
\end{eqnarray*}
We have that global sections of $\left( \widetilde \omega^{\Gamma}_0 \right)^{\vee}$ are sections of $\Of_{\Pro^1_{\C}}(1)$ on extremal components and $\Of_{\Pro^1_{\C}}(-1)$ on the central component which agree on nodes. Since $H^0(\Pro^1_\C,\Of_{\Pro^1_{\C}}(-1))= 0$,  global sections of $\left( \widetilde \omega^{\Gamma}_0 \right)^{\vee}$ are sections of $\Of_{\Pro^1_{\C}}(1)$ on the three extremal components which are zero on nodes. On each extremal component we fix coordinates $[X'_0,X'_1]$, $[X''_0,X''_1]$ and $[X'''_0, X'''_1]$. Sections on $\Of_{\Pro^1_{\C}}(1)$ are linear forms
$$ 
a_1 X'_0 + b_1 X'_1 \quad a_2 X''_0 + b_2 X''_1 \quad a_3 X'''_0 + b_3 X'''_1
$$
which are zero on nodes if and only if $a_1=a_2=a_3=0$. We have that only $(\Gm)^{\Delta_1}$ does not acts trivially on $H^0(C_0,\left( \widetilde \omega^{\Gamma}_0 \right)^{\vee})$ and the action is
\begin{eqnarray*}
\az: (\Gm \times \Gm \times \Gm) \times H^0(C_0,\left( \widetilde \omega^{\Gamma}_0 \right)^{\vee}) &\to& H^0(C_0,\left( \widetilde \omega^{\Gamma}_0 \right)^{\vee})\\
(\lambda_1,\lambda_2, \lambda_3), (b_1, b_2, b_3 ) &\mapsto& ( \lambda_1 b_1, \lambda_2 b_2, \lambda_3 b_3)
\end{eqnarray*}
so the Chern character of $\widetilde \Pi_* \omega^{\vee}_0$ is
$$
e^{t_1} + e^{t_2} + e^{t_3} = 3 + \sum^{\infty}_{m=1} \frac{t^m_1+ t^m_2 + t^m_3}{m!}
$$
as we claimed.
\end{dem}

\section{The Chow ring of $\M_0^{\leq 3}$}

In this Section we calculate $A^*(\M^{\leq 3}_0) \rat$.

\subsection{The open substack $\M _0 ^0$}

\medskip

\framebox{$\Gamma:= \bullet$}

\medskip

From the equivalence (\ref{equiv})
\begin{equation}
\M_0 ^{\Gamma} \simeq \B \aut(C).
\end{equation}

we have that the stack $\M_0 ^0$ is the classifying space of $\Pro Gl_2$. Owing to the fact that $\Pro Gl_2 \cong SO_3$ and following \cite{Pa} we have
$$
A^*(\M_0 ^0) \rat  \cong \Q[c_2(\sla _2)].
$$
Since (see Proposition \ref{mumzero}) $c_2(\sla _2)=(1/2)\km_2$, we can write
\begin{propos}
$$
A^*(\M_0 ^0) \rat =\Q[\km_2].
$$
\end{propos}

\subsection{The first stratum}

\framebox{$\Gamma:=\; \small{\unoa}$} 

\medskip

We order the two components. The automorphism group is $\Cg_2 \ltimes (E \times E)$, (where $\Cg_2$ is the order two multiplicative group) and the action of its generator $\tau$ over  $E \times E$ exchanges the components.

Then the induced action of $\tau$ on the ring $A^*_{E \times E} \rat \cong \Q [t_1, t_2]$ (see Proposition \ref{pr}) exchanges the first Chern classes $t_1$ and $t_2$. The invariant polynomials are the symmetric ones which are algebrically generated by: $\{ (t_1 + t_2)/2, (t_1^2  + t_2 ^2)/2 \}$. By recalling the description (\ref{kmum}) of Mumford classes, we have $A^*(\M _0 ^1) \rat = \Q[\km_1, \km_2]$.
Let us consider the two inclusions $i$ and $j$
(respectively closed and open immersions) and the \'etale covering $\phi$
\begin{equation}
\xymatrix{
\widetilde \M_0^1 \ar[r]^{\phi}&\M_0 ^1 \ar [dr]^i \\
&\M_0 ^0 \ar [r]_j & \M_0 ^{\leq 1}
}
\end{equation}
we obtain the following exact sequence
$$
A^*(\M_0 ^1) \rat \xrightarrow{i_*} A^*(\M_0 ^{\leq 1}) \rat \xrightarrow{j^*} A^*(\M_0 ^0) \rat \xrightarrow{} 0
$$
for what we have seen we have:
$$
\Q[\km_1, \km_2] \xrightarrow{i_*} A^*(\M_0 ^{\leq 1}) \rat \xrightarrow{j^*} \Q[\km_2] \xrightarrow{} 0.
$$
Now with reference to the paragraph (\ref{normal}) we have that the first Chern class of the normal bundle $N_{\M_0^1} ( \M_0^{\leq 1})$ is
$$
i^*i_* [\M_0^1]=\phi_*\frac{1}{2}(t_1+t_2)=-\km_1
$$
Since $A^*{\M^1_0}$ is an integral domain we can apply Lemma (\ref{incollamento}) and obtain the ring isomomorphism
$$
A^*(\M_0 ^{\leq 1}) \rat \cong \Q[\km_1,\km_2] \times_{\Q[\km_2]} \Q[\km_2] \cong \Q[\km_1, \km_2].
$$
where the map $q: \Q[\km_2] \to \Q[\km_1,\km_2]/(\km_1)=\Q[\km_2]$ tautologically sends $\km_2$ into $\km_2$.

So we have
\begin{propos}
$$
A^*(\M_0 ^{\leq 1}) \rat \cong \Q[\km_1, \km_2].
$$
\end{propos}
\subsection{The second stratum}

\framebox{$\Gamma:=\;\small{\duea}$} We order the two components with one node.

In this case the group of automorphism of the fiber is
$$
\text{Aut}(C^{\Gamma}) \cong \Cg_2 \ltimes (\Gm \times E \times E)=:\Cg_2 \ltimes H,
$$
where the action of $\tau$ sends an element $g \in \Gm$ into $g^{-1}$ and exchange the components isomorphic to $E$.

We can identify $A^*(\B (\Gm)^3)$ with $A^*(\widetilde{\M}_0^{\Gamma})$ and $A^*(\B \aut(\Gamma) \ltimes (\Gm)^3)$ with $A^*(\M_0^{\Gamma})$.  Following the notation of Section \ref{stratclass} set
\begin{eqnarray*}
t_1=\psi^1(\infty,1) \quad t_2=\psi^1(\infty,2) \quad r=\psi^2(\infty)&
\end{eqnarray*}
the action induced by $\tau$ on these classes is $\tau(r, t_1, t_2) = (-r, t_2, t_1)$.
With reference to the map
$$
\B (\Gm)^3 \xrightarrow{\phi} \B \aut(\Gamma) \ltimes (\Gm)^3
$$
we recall that $\phi^*$ is an isomorphism between $A^*(\M_0^1)\rat$ and $A^*(\B (\Gm)^3)^{\Cg_2}$ (see Proposition \ref{isophi}).

We can describe $A^*(\B H) \rat=\Q [r, t_1, t_2]$ as the polynomial ring in $r$ with coefficients in $\Q[t_1,t_2]$, so we write a polynomial  $P(r,t_1,t_2)$ as $\sum _{i=0} ^k r^i P_i(t_1, t_2)$.

The polynomial $P$ is invariant for the action of $\tau$ if and only if the coefficients of the powers of $r$  in $P(r,t_1,t_2)$ are equal to those of the polynomial $P(-r,t_2,t_1)$.

That is to say that $P_i$ with even index are invariant for the exchange of $t_1$ and $t_2$, while those with odd index are anti-invariant. An anti-invariant polynomial $Q$ is such that $Q(t_1,t_2) + Q(t_2, t_1)=0$
and consequently it is the product of $(t_1 - t_2)$ by an invariant polynomial. It is furthermore straightforward verifying that any such polynomial is invariant for the action of $\tau$.

So an algebraic system of generators for $(A^*_{(\Gm)^3}\rat)^{\Cg _2}$ is given by
$$
u_1:=t_1 + t_2 \quad u_2:=t_1 ^ 2 +  t_2 ^2, \quad u_3:=r(t_1-t_2), \quad u_4:=r^2.
$$
We know that (see Sections \ref{mumcla} and \ref{stratclass})
\begin{eqnarray*}
&&\phi^*\km_1=- u_1\\
&&\phi^*\km_2=- u_2\\
&&\phi^*(\gamma_2) = (t_1-r)(t_2+r)=\frac{1}{2}(u_1^2-u_2) +u_3 - u_4
\end{eqnarray*}
where $\gamma_2=c_2(\Nf_{\M_0^2/\M_0^{\leq 2}})$, and there exists a class $x \in A^2\M_0^2\rat$ such that $u_3=\pi^*x$.
\begin{claim}
The ideal of relations is generated on $\Q [\km_1, \km_2, \gamma_2, x]$ by the polynomial
\begin{eqnarray} \label{p2a}
(2x + (2\km_2+\km_1^2))^2-(2\km_2 + \km_1^2)(4\gamma_2 - \km_1^2)=0.
\end{eqnarray}
\end{claim}

\begin{dem}
From direct computation we have that relation (\ref{p2a}) holds and the polynomial is irreducible. On the other hand let us consider the map $f: \A^3_{\Q} \to \A^4_{\Q}$ defined as $(t_1,t_2,r) \mapsto (u_1, u_2, u_3, u_4)$.
If the generic fiber of $f$ is finite then $f(\A^3_{\Q})$ is an hypersurface in $\A^4_{\Q}$ and we have done.
Now for semicontinuity it is sufficient to show that a fiber is finite.
Let us consider the fiber on 0. We have that $u_1, u_2, u_3, u_4$ are simultaneously zero iff $t_1=t_2= r=0$.
\end{dem}

{\bf NOTE:} In the following we do not explicit the argument above.

\bigskip

Now set $\eta:=2x + (2k_2 + k_1^2)$,
we have that $A^*(\M _0 ^2) \rat$ is isomorphic to the graded ring $\Q [\km_1,\km_2,\gamma_2, \eta]/I$
where the ideal $I$ is generated by the polynomial $\eta^2-(2\km_2 + \km_1^2)(4\gamma_2 - \km_1^2)$.
Since $\phi^* \phi_*$ is multiplication by two, we also have the following relation
$$
\eta=\phi_*\left( \frac{1}{2}(t_1-t_2)(2r-t_1+t_2)\right).
$$

Let us consider the cartesian diagram
\begin{equation*}
\xymatrix{
A^*(\M_0^{\leq 2}) \rat  \cart \ar[r]^{j^*} \ar[d]_{i^*} &\Q[\km_1, \km_2] \ar[d]^{q}\\
\Q[\km_1, \km_2, \gamma_2, \eta]/I \ar[r]^{p} & \Q[\km_1, \km_2, \eta]/\overline{I}
}
\end{equation*}
where $\overline{I}$ is the ideal generated by $\eta^2 +(2\km_2 +\km_1^2)\km_1^2$.

The map $q$ is injective so $i^*$ is injective too.

We set in $A^*(\M_0^{\leq 2}) \rat$ the classes $\gamma_2:=i_*1$ and $q:=i_*\eta$,
the ring we want (from injectivity of $i^*$) is isomorphic  to the subring of $A^*(\M_0^2)\rat$ generated by $\km_1, \km_2, \gamma_2, \gamma_2 \eta$ so we have
\begin{propos}
$$
A^*(\M_0^{\leq 2}) \rat = \frac{\Q [\km_1, \km_2, \gamma_2, q]}{(q^2 +\gamma_2 ^2(2\km_2 + \km_1^2)(\km_1^2-4\gamma_2))}
$$
\end{propos}
\subsection {The third stratum}

The third stratum splits into two components.

\bigskip

{\bf The first component}
$\;$

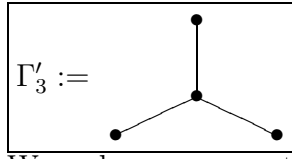
\framebox{$\Gamma'_3:=\;\vcenter{\small{\treb}}$}

We order components in $\Delta_1$. Let us note that the component corresponding to the central vertex has three points fixed by the other three components, consequently, given a permutation of the external vertices, there is an unique automorphism related to the central vertex.
 
The group $\text{Aut}(C^{\Gamma'_3})$ is therefore isomorphic to $S_3 \ltimes (E^3)$.

As usual from proposition (\ref{pr}) it follows that
$$
A^*_{E^3} \rat \cong \Q [w_1, w_2, w_3],
$$
on which $S_3$ acts by permuting the three classes
\begin{eqnarray*}
w_1:=\psi^1_{\infty,1} \quad w_2:=\psi^1_{\infty,2} \quad w_3:=\psi^1_{\infty,3}
\end{eqnarray*}
So we have
$$
\begin{array}{l}
\phi^*\km_1=-(w_1 + w_2 + w_3),\\
\phi^*\km_2:=-(w_1^2 + w_2^2 + w_3^2),\\
\phi^*\km_3:=-(w_1^3 + w_2^3 + w_3^3);
\end{array}
$$
conesequently
$$
A^*(\M_0 ^{\Gamma'_3}) \rat \cong \Q [\km_1, \km_2, \km_3].
$$

As we have seen in Section \ref{restrcla} we fix the following notation
\begin{equation*}
\xymatrix@=1.3pc{
\widetilde \M_0^{\Gamma'_3}   \ar@/_2pc/[dddrr]_\phi \ar[rr]^-{f}& &\Psi(\Gamma_2, \Gamma'_3)\ar[rr]^{pr_2} \ar [ddd]^{pr_1}& & \left( \widetilde{\M_0^{\Gamma_2}} \right) ^{\leq 2} \ar[ddd]^{\Pi}\\
& &\\
& & & \\
& & \M_0^{\Gamma'_3} \ar[rr]^{in}& &\M_0^{\leq 3}
}
\end{equation*}
where $f$ is the union of all $f_{\alpha}$.

First of all let us notice that the class $\phi^*\gamma_3':=\pi^*c_3(\Nf_{i})=w_1 w_2 w_3$
depends on Mumford classes in the following way $6\gamma'_3=-(\km_1^3 - 3 \km_1 \km_2 + 2 \km_3)$
so we can write
$$
A^*(\M_0^{\Gamma'_3}) \rat = \Q[\km_1, \km_2, \gamma'_3].
$$
The restriction of $\gamma_2$ to $A^*(\M_0^{\Gamma'_3})$is
\begin{eqnarray*}
\gamma_2:=pr_{1*}\frac{1}{2}c_2(pr_1^*(\Nf_{in})/\Nf_{pr_2})&=&\phi_*f^*\frac{1}{2}c_2(pr_1^*(\Nf_{in})/\Nf_{pr_2})\\
&=&\frac{1}{2}\pi_*(w_1 w_3)
\end{eqnarray*}
from which $\phi^*\gamma_2=w_1 w_3 + w_1 w_2 + w_2 w_3$ consequently, by writing $2\gamma_2 = \km_1^2 - \km_2$ we have
$$
A^*(\M_0^{\Gamma'_3}) \rat = \Q[\km_1, \gamma_2, \gamma'_3].
$$
In order to restrict the class $q$ let us notice that we can write
\begin{eqnarray*}
f^*r = 0 \quad f^*t_1= w_1 \quad f^*t_2 = w_3
\end{eqnarray*}
from which we have
$$
f^* \left( \frac{1}{2}(t_1 - t_2)(2t - t_1 + t_2) \right)= - \frac{1}{2}(w_1 - w_3)^2
$$
and so
\begin{eqnarray*}
\phi^*q&=&\phi^*\phi_* \left( - \frac{1}{2}(w_1 - w_3)^2 w_1w_3 \right)\\
&=& - ( (w_1 - w_3)^2 w_1 w_3 + (w_1 - w_2)^2 w_1 w_2 + (w_2 - w_3)^2 w_2 w_3)
\end{eqnarray*}
we can therefore write $q= - \gamma_2(\km_1^2 - 4 \gamma_2) + 3 \gamma'_3 \km_1$.
With reference to the inclusions
$$
\M_0^{\Gamma'_3} \xrightarrow{i} \M_0^{\leq 2} \cup \M_0^{\Gamma'_3} \xleftarrow{j} \M_0^{\leq 2}
$$
we have the fiber square:
\begin{equation*}
\xymatrix{
A^*(\M_0^{\leq 2} \cup \M_0^{\Gamma'_3}) \rat  \cart \ar[r]^{j^*} \ar[d]_{i^*} &\Q[\km_1, \km_2, \gamma_2, q]/I \ar[d]^{\varphi}\\
\Q[\km_1, \gamma_2, \gamma'_3] \ar[r]^{p} & \Q[\km_1, \gamma_2]
}
\end{equation*}
where  $I$ is the ideal generated by the polynomial $q^2 - \gamma_2^2(2\km_2 - \km_1^2)(\km_1^2 - 4\gamma_2)$ and the map $\varphi$ is surjective and such that
$$
\ker \varphi = (q + \gamma_2(\km_1^2 - 4 \gamma_2), \quad \km_2 + 2\gamma_2 - \km_1^2).
$$ 
Now let us observe that from Lemma (\ref{quoz}) the ring in question is isomorphic to
$$
A/(0,I):=\frac{\Q[\km_1, \gamma_2, \gamma'_3] \times_{\Q[\km_1, \gamma_2]} \Q[\km_1, \km_2, \gamma_2, q] }{(0,I)}.
$$
Set, with abuse of notation
\begin{eqnarray*}
\km_1:=(\km_1,\km_1) \quad & \gamma_2:=(\gamma_2,\gamma_2) \quad & \gamma'_3:=(\gamma'_3, 0)\\
\km_2:=(\km_1^2 - 2\gamma_2, \km_2) \quad & q:=(-\gamma_2(\km_1^2 - 4\gamma_2),q)&
\end{eqnarray*}
Straightforward arguments lead us to state the following
\begin{lemma}
The classes $\km_1,\gamma_2, \gamma'_3, \km_2, q$ generate the ring $A$.
\end{lemma}
Let us compute the ideal of relations. Let $T(\km_1, \gamma_2, \gamma'_3, \km_2, q)$ be a polynomial in $\Q[\km_1, \gamma_2, \gamma'_3, \km_2, q]$, it is zero in $A$ iff $T(\km_1, \gamma_2, \gamma'_3, \km_1^2 - 2\gamma_2, -\gamma_2(\km_1^2 - 4\gamma_2))$ and $T(\km_1, \gamma_2, 0, \km_2, q)$ are respectively zero in $\Q[\km_1, \gamma_2, \gamma'_3]$ and $\Q[\km_1, \gamma_2, \km_2, q]$: in particular this imply that $T$ is in the ideal of $\gamma'_3$. Consequently the polynomial $T=:\gamma'_3 \widehat{T}$ is zero in $A$ iff $\widehat{T}(\km_1, \gamma_2, \gamma'_3, \km_1^2 - 2\gamma_2, -\gamma_2(\km_1^2 - 4\gamma_2))$ is zero in  $\Q[\km_1, \gamma_2, \gamma'_3]$. The ideal of relations in $A$ is so generated by $\gamma'_3(-\km_2 + 2\gamma_2 - \km_1^2)$ and $\gamma'_3(q + \gamma_2(\km_1^2- 4\gamma_2))$.

Finally let us notice that the ideal $(0,I)$ is generated in $A$ by the polynomial $ q^2 + \gamma_2^2(2\km_2 + \km_1^2)(\km_1^2 - 4\gamma_2)$. So we can conclude that
$$
A^*(\M_0^{\leq 2} \cup \M_0^{\Gamma'_3}) \rat = \Q[\km_1, \km_2, \gamma_2, \gamma'_3, q]/J,
$$
where $J$ is the ideal generated by the polynomials
\begin{eqnarray*}
&& q^2 + \gamma_2^2(2\km_2 + \km_1^2)(\km_1^2 - 4\gamma_2)\\
&& \gamma'_3(- \km_2 + 2\gamma_2 - \km_1^2)\\
&&\gamma'_3(q + \gamma_2(\km_1^2- 4\gamma_2))
\end{eqnarray*}

{\bf The second component}

\framebox{$\Gamma''_3:=\;\vcenter{\small{\trea}}$}

The group of $\text{Aut}(C^{\Gamma'_3})$ is $\Cg_2 \ltimes (E \times E \times \Gm \times \Gm)$. The action of $\tau$ on this group exchange  simultaneously the components related to $\Gm$ and those related to $E$.

We have the isomorphism
$$
A^*_{ E \times E \times \Gm \times \Gm} \rat \cong \Q [v_1, \dots, v_4]
$$ 
where
\begin{eqnarray*}
v_1=\psi^1_{\infty,1} \quad v_2=\psi^1_{\infty, 2} \quad v_3=\psi^2_{\infty, 1} \quad v_4&=&\psi^2_{\infty, 2}
\end{eqnarray*}
by gluing curves such that the two central components corresponds in the point at infinity.

It follows that the action induced by $\tau$ is $\tau(v_1,v_2,v_3,v_4)=(v_2, v_1,v_4,v_3)$.
Since $\Cg_2$ has order 2, the invariant polynomials are algebraically generated by the invariant polynomials of degree at most two (see. Theorem 7.5 \cite{CLO}). It is easy to see that a basis for the linear ones is given by $u_1:= v_1 + v_2, \; u_2:=v_3 + v_4$. For the vector subspace of invariant polynomials of degree two, we can compute a linear basis by using Reynolds' operator
$$
 \begin{array}{ll}
u_3:=v_1 ^2 +  v_2 ^ 2, & u_6:=v_1 v_3 + v_2 v_4,\\ 
u_4:=v_3 ^2 + v_4 ^2, & u_7:=v_1 v_2,\\
u_5:=v_1 v_4 + v_2 v_3, & u_8:=v_3 v_4
\end{array}
$$
now we note that
$$
u_6= u_1 u_2 - u_5,\quad u_7= (u_1^2 - u_3)/2,\quad u_8= (u_2^2 - u_4)/2
$$
Consequently we can write
$$
A^* \M_0^{\Gamma'_3} \rat \cong \Q [u_1, \dots, u_5]/I,
$$
where $I$ is the ideal generated by the polynomial
\begin{equation} \label{polgammauno}
2u_3 u_4 + 2u_1 u_2 u_5 -u_2 ^2 u_3 -u_1^2 u_4 - 2u_5^2
\end{equation}
With reference to the degree two covering  $\widetilde{\M}_0^{\Gamma'_3} \xrightarrow{\phi} \M_0^{\Gamma'_3}$
we have: $- u_1=\phi^*\km_1 \; \; -u_3=\phi^*\km_2$.
In order to compute the restriction of the closure of the classes $\gamma_2$ and $q$ of $\M_0^{\leq 2}$ let us fix the notation of the following diagram
\begin{equation*}
\xymatrix@=2.0pc{
& &\Ms_{0,1}^1 \times \Ms_{0,2}^0 \times \Ms_{0,1}^0 \ar[drr]^{pr_2^{'}} \\ \widetilde{\M}_0^{\Gamma'_3} \ar[urr]^{f_1} \ar[drr]^{f_3} \ar@/_2pc/[dddrr]_\phi \ar[rr]^-{f_2}& &(\Ms_{0,1}^0 \times \Ms_{0,2}^1 \times \Ms_{0,1}^0)^I \ar[rr]^{pr_2^{''}} & & \left( \widetilde{\M}_0^{\Gamma_2} \right)^{\leq 2} \ar[ddd]^{\Pi}\\
& &\Ms_{0,1}^0 \times \Ms_{0,2}^0 \times \Ms_{0,1}^1 \ar[dd]^{pr_1}\ar[urr]^{pr_2^{'''}}\\
& & & \\
& & \M_0^{\Gamma''_3} \ar[rr]^{i}& &\M_0^{\leq 3}
}
\end{equation*}
where $(\Ms_{0,1}^0 \times \Ms_{0,2}^1 \times \Ms_{0,1}^0)^I $ is the component where the marked points of the central curve (which is singular) are on different components.

As $\phi^*$ is an isomorphism to the algebra of polynomials which are invariants for the action of
$\Cg_2$, let us choose in $A^*(\M_0^{\Gamma''_3})\rat$ classes $\rho, \lambda, \mu$ such that
$u_2=\phi^*\rho \; \;,
u_4=\phi^*\lambda \; \;,
u_5=\phi^*\mu$.

First of all let us compute the restriction of the closure of $\gamma_2 \in A^*(\M_0^{\leq 2}) \rat$, the polynomial in the classes $\psi$ is $\frac{1}{2}$, we have
$c_3(\phi^*(\Nf_{i})) =(v_1-v_3)(v_3+v_4)(v_2-v_4)$, $c_1(f_1^*(\Nf_{pr_2^{'}}))=(v_1-v_3)$, $c_1(f_2^*(\Nf_{pr_2^{''}}))=(v_3+v_4)$ and $c_1(f_3^*(\Nf_{pr_2^{'''}}))=(v_2-v_4)$
from which we obtain the following relations
\begin{eqnarray*}
&\phi^*\gamma'_3&=(v_1-v_3)(v_3+v_4)(v_2-v_4)\\
&&=\phi^*\left(\rho \left(\frac{1}{2}\left(\km_1^2 + \rho^2 - \km_2 - \lambda \right) -\mu \right) \right)\\
&\phi^*\gamma_2&=\left( (v_3+v_4)(v_2-v_4) + (v_1-v_3)(v_2-v_4) + (v_1-v_3)(v_3+v_4) \right)\\
&&=- \rho \km_1 - \mu + \frac{1}{2}(\km_1^2 - \rho^2 + \km_2 - \lambda)
\end{eqnarray*}
In order to have $\gamma_2$ among the generators, set
$$
\lambda= -2 \rho \km_1 - 2\mu + \km_1^2 - \rho^2 + \km_2 - 2 \gamma_2;
$$
the equation (\ref{polgammauno}) becomes
\begin{equation} \label{polbis}
K: \; \; \sigma^2-(2\km_2 + \km_1^2)((-\km_1 + 3\rho)(\km_1 + \rho) -4\gamma_2)
\end{equation}
where we have set $\sigma=2\mu - 2\km_2 - \km_1^2 + \km_1 \rho$.
Then we can write the ring $A^*(\M_0^{\Gamma''_3})\rat$ as $\Q[\km_1, \rho, \km_2, \gamma_2, \sigma]/K$.
In the new basis we have $\gamma'_3=\rho( \rho^2 - \rho \km_1 + \gamma_2)$.
Let us restricts the closure of $q \in A^*(\M_0^{\leq 2})\rat$ to $\M_0^{\Gamma_1}$. We recall that the related polynomial in classes $\psi$ of $\widetilde{\M}_0^{\Gamma_2}$ is $\frac{1}{2}(t_1-t_2)(2r - t_1+ t_2)$.
On $\Ms_{0,1}^1 \times \Ms_{0,2}^0 \times \Ms_{0,1}^0$ we have $
f_1^*t_1=v_3 \quad f_1^*t_2 = v_2 \quad f_1^*r =-v_4,$
on $(\Ms_{0,1}^0 \times \Ms_{0,2}^1 \times \Ms_{0,1}^0)^I$ we have $f_2^*t_1= v_1\quad f_2^*t_2 = v_2 \quad f_2^*r =\frac{\psi^2_\infty -\psi^2_0}{2}=\frac{v_3 - v_4}{2},$
and in the end on $\Ms_{0,1}^0 \times \Ms_{0,2}^0 \times \Ms_{0,1}^1$ we have $f_3^*t_1=v_1 \quad f_3^*t_2 =v_4 \quad f_3^*r =v_3,$
consequently the polynomials related to the three components of $\Psi(\Gamma_2, \Gamma''_3)$ are
$P_1 = \frac{1}{2}(v_3-v_2)(v_2 - v_3 - 2v_4)$, $P_2= \frac{1}{2}(v_1 - v_2) \left( v_4 - v_3 - v_1 + v_2 \right)$ and $P_3 = \frac{1}{2}(v_4 - v_1)(v_1 - 2v_3 - v_4)$
from which
\begin{eqnarray*} 
\phi^*q&=& 2\sum_{\alpha=1}^3 (P_\alpha N_\alpha)
= \phi^*((3\rho + \km_1)\gamma''_3 + (\rho^2 - \gamma_2)\sigma)
\end{eqnarray*}
this means that the image of $q$ in $B:=A^*(\M_0^{\Gamma''_3}) \rat / {\gamma''_3}$ is $(\rho^2 - \gamma_2)\sigma$.
Further let us notice that the image of $\gamma'_3$ in $B$ is 0.

In order to compute $A^*(\M_0^{\leq 3}) \rat$, let us consider its isomorphism with the fibered product
\begin{equation*}
\xymatrix{
 A^*(\M_0^{\leq 3}) \rat \cart \ar[r]^{j^*} \ar[d]_{i^*} & \Q [\km_1, \km_2, \gamma_2, \gamma'_3, q]/J \ar[d]^{\varphi}\\
\Q[\km_1, \rho, \km_2, \gamma_2, \sigma]/K \ar[r]^{p} & B
}
\end{equation*}
With reference to the fiber square
\begin{equation*}
\xymatrix{
 A \cart \ar[r]^{j^*} \ar[d]_{i^*} & \Q [\km_1, \km_2, \gamma_2, \gamma'_3, q] \ar[d]^{\varphi}\\
\Q[\km_1, \rho, \km_2, \gamma_2, \sigma]/K \ar[r]^{p} & B
}
\end{equation*}
we have that the ring in question is isomorphic to the quotient $A/(0,J)$.

Now we look for generators of the ring $A$.

Again a straightforward argument leads us to state
\begin{lemma}
The following elements of $A$
$$
\begin{array}{lll}
\km_1:=(\km_1, \km_1); & \km_2:=(\km_2,\km_2); & \gamma_2:=(\gamma_2);\\
\gamma'_3:=(0,\gamma'_3); & q:=((\rho^2 - \gamma_2)\sigma,q); & \gamma''_3:=(\gamma''_3,0);\\
r:=(\gamma''_3 \rho,0); & s:=(\gamma''_3 \sigma, 0) & t:=(\gamma''_3 \rho^2, 0);\\
& u:=(\gamma''_3 \rho \sigma,0)&
\end{array}
$$
are generators of the ring.
\end{lemma}

Now let us call $\widetilde A$ the ring $\Q[\km_1, \km_2, \gamma_2, \gamma'_3, \gamma''_3, q, r,s,t,u]$. We have by fact defined a surjective homomorphism $a: \widetilde A \to A$; we call again $j^*$ and $i^*$ their composition with $a$, we have $\ker a = \ker i^* \cap \ker j^*$.
Now, for what we've seen, we have:
\begin{eqnarray*}
\ker i^*&=& \left( 
	\begin{array}{l} 
	\gamma'_3;\\ 
	r^2 - \gamma''_3t; \\
	rs - \gamma''_3 u;\\
	s^2 - (\gamma''_3)^2(2\km_2 + \km_1^2)((-\km_1 + 2\rho)(\km_1 + \rho) + 4\gamma_2)
	\end{array}
	\right)\\
\ker j^* &=& (\gamma''_3, r,s,t,u)
\end{eqnarray*}
and so:
$$
\ker a= \left(
\begin{array}{l}
\gamma'_3 \gamma''_3; \quad \gamma'_3 r; \quad \gamma'_3 s; \quad \gamma'_3 t; \quad \gamma'_3 u;\\
r^2 - \gamma''_3t; \quad rs - \gamma''_3 u;\\
s^2 - (\gamma''_3)^2(2\km_2 + \km_1^2)((-\km_1 + 2\rho)(\km_1 + \rho) + 4\gamma_2)
\end{array}
\right)
$$

We make the quotient of $A=A'/ \ker a$ by the ideal $(0,J)$ and we obtain the following:
\begin{theorem}\label{finthr}
The ring $A^*(\M_0 ^{\leq 3}) \rat$ is: $\Q[\km_1, \km_2, \gamma_2, \gamma'_3, \gamma''_3, q, r,s,t,u]/L,$
where $L$ is the ideal generated by the polynomials
$$
\begin{array}{l}
\gamma'_3(- \km_2 + 2\gamma_2 - \km_1^2); \quad \gamma'_3(q + \gamma_2(\km_1^2 - 4 \gamma_2));\\
\gamma'_3 \gamma''_3; \quad \gamma'_3 r; \quad \gamma'_3 s; \quad \gamma'_3 t; \quad \gamma'_3 u;\\
r^2 - \gamma''_3t; \quad rs - \gamma''_3 u;\\
q^2 + (\gamma_2)^2(2 \km_2 + \km_1^2)(\km_1^2 - 4 \gamma_2);\\
s^2 - (2\km_2 + \km_1^2)((- \km_1\gamma''_3 + 2r)(\km_1\gamma''_3 + r) + 4\gamma_2(\gamma''_3) ^2)
\end{array}
$$

\end{theorem}


\begin{thebibliography}{99}

\bibitem[\textbf{Art}]{Artd} M.Artin: Versal Deformation and Algebraic Stacks. Invent. Math. \textbf{27}, 165-189 (1974).

\bibitem[\textbf{CLO}]{CLO} D.Cox, J.Little, D.O'Shea: \emph{Ideals, Varieties, and Algorithms} Springer-Verlag: New York, 1998.

\bibitem[\textbf{Ed-Gr}]{EGRR} D.Edidin, G.Graham: Equivariant intersection theory; Inv. Math. \textbf{131}, 595-634 (1998).

\bibitem[\textbf{EHKV}]{EHKV} D.Edidin, B.Hassett, A.Kresh, A.Vistoli: Brauer groups and quotient stacks;  Amer. J. Math.  \textbf{123}  no. 4, 761-777 (2001).

\bibitem[\textbf{EGA IV}]{egaq} A.Grothendieck and J.Dieudonn\'e: \'Etude locale des sh\'emas et des morphismes de sch\'emas; Publ. Math. I.H.E.S. \textbf{20} (1964), \textbf{24} (1964), \textbf{28} (1965), \textbf{32} (1966).

\bibitem[\textbf{Fab}]{Fab} C.Faber: Chow rings of moduli spaces of curves; Ann. of Math (2) \textbf{132} 331-449 (1990).

\bibitem[\textbf{Ful}]{ful} W.Fulton: \emph{Intersection theory} Springer-Verlag (Berlin), 1998.

\bibitem[\textbf{Fulg}]{fulg} D.Fulghesu, PhD Thesis, Scuola Normale Superiore, Pisa 2005.

\bibitem[\textbf{Gil}]{Gil} H.Gillet: Intersection theory on algebraic stacks and $Q$-varieties; J. Pure Appl. Algebra \textbf{34}, 193-240 (1984).

\bibitem[\textbf{Ha-Mo}]{HM} J.Harris, I.Morrison: \emph{Moduli of curves} Springer-Verlag (New York), 1998.

\bibitem[\textbf{Hir}]{Hir} H.Hironaka: An example of a non-Kahlerian deformation, Ann. of Math. \textbf{75}, 190- (1962).

\bibitem[\textbf{Hrt}]{hrt} R.Hartshorne: \emph{Algebraic Geometry} Springer-Verlag New York, 1997.

\bibitem[\textbf{Knu}]{Kn} D.Knutson: \emph{Algebraic Spaces, Lecture Notes in Math.}, vol. 203, Springer-Verlag (Berlin), 1971.

\bibitem[\textbf{Kre}]{kre} A.Kresch: Cycle groups for Artin stacks; Inv. Math. \textbf{138}, 495-536 (1999).

\bibitem[\textbf{Lang}]{lan} S.Lang: \emph{Algebra} Springer-Verlag, New York (2002).

\bibitem[\textbf{Mum1}]{mumlect} D.Mumford: \emph{Lectures on curves on an algebraic surface} Princeton University Press: Princeton, 1966.

\bibitem[\textbf{Mum2}]{mumenum} D.Mumford: Towards an Enumerative Geometry of the Moduli Space of Curves; In \emph{Arithmetic and geometry, II}, \textbf{36} Progress in Math., 271-326: Boston, 1983.

\bibitem[\textbf{Pan}]{Pa} R.Pandharipande: Equivariant Chow rings of $O(k), SO(2k+1)$ and $SO(4)$; J. Reine Angew. Math. \textbf{496}, 131-148 (1998).

\bibitem[\textbf{Shf}]{Shf} I.R.Shafarevich: \emph{Basic Algebraic Geometry 2} Springer-Verlag (Berlin) 1994.

\bibitem[\textbf{Tot}]{tot} B.Totaro: The Chow ring of a classifying space; Proc. Sympos. Pure Math., \textbf{67}, Amer. Math. Soc., Providence, RI, (1999).

\bibitem[\textbf{Vez}]{Ve} G.Vezzosi: On the Chow ring of the classifying stack of ${\rm PGL}_{3,\C}$; J. Reine Angew. Math. \textbf{523} 1-54 (2000).

\bibitem[\textbf{Ve-Vi}]{VV} G.Vezzosi, A.Vistoli: Higher algebraic $K$-theory for actions of diagonalizable groups; Invent. math. \textbf{153}, 1-44 (2003).

\bibitem[\textbf{Vis}]{Vis} A.Vistoli: Intersection theory on algebraic stacks and their moduli spaces; Invent. Math. \textbf{97}, 613-670 (1989).

\end{thebibliography}
\end{document}